\documentclass[MSNbibl,number,citesort,seceqn,dvips]{arxbj}
\usepackage{subenv}

\aid{0}
\volume{19}
\issue{5A}
\pubyear{2013}
\firstpage{1612}
\lastpage{1636}
\doi{10.3150/12-BEJ423} 

\makeatletter
\newtheorem{theorem}{Theorem}[section]
\newremark{remark}{Remark}
\newremark{example}[theorem]{Example}
\newtheorem{corollary}[theorem]{Corollary}
\makeatother

\begin{document}
\begin{frontmatter}

\title{Asymptotic mean stationarity and absolute continuity of point
process distributions}
\runtitle{Asymptotic mean stationarity}

\begin{aug}
\author{\fnms{Gert} \snm{Nieuwenhuis}\corref{}\ead[label=e1]{g.nieuwenhuis@uvt.nl}}
\runauthor{G. Nieuwenhuis} 
\address{Department of Econometrics \& OR,
Tilburg School of Economics and Management (TiSEM),
Tilburg University,
P.O. Box 90153, 5000 LE Tilburg, The Netherlands.\\ \printead{e1}}
\end{aug}

\received{\smonth{4} \syear{2011}}
\revised{\smonth{12} \syear{2011}}

%
\begin{abstract}
This paper relates -- for point processes $\Phi$ on $\mathbb{R}$ --
two types of asymptotic mean stationarity (AMS) properties and several
absolute continuity results for the common probability measures emerging
from point process theory. It is proven that $\Phi$ is AMS under the
time-shifts if and only if it is AMS under the event-shifts. The
consequences for the accompanying two types of ergodic theorem are
considered. Furthermore, the AMS properties are equivalent or closely
related to several absolute continuity results. Thus, the class of AMS
point processes is characterized in several ways. Many results from
stationary point process theory are generalized for AMS point processes.
To obtain these results, we first use Campbell's equation to rewrite the
well-known Palm relationship for general nonstationary point processes
into expressions which resemble results from stationary point process
theory.
\end{abstract}

%

\end{frontmatter}

\section{Introduction}\label{sec1}

Point process theory on $\mathbb{R}$ utilizes two types of shifts
(event-shifts and time-shifts) and several closely related probability
measures. Each type of shifts brings its own ergodic theorem and -- by
taking expectations -- its own concept `asymptotic mean stationarity'
(AMS). This paper develops a theory of AMS point processes on
$\mathbb{R}$. The basic result is that the two types of AMS are
equivalent, thus extending a classical result of Kaplan \cite{K55} about
equivalence of event-stationarity and time-stationarity. Furthermore,
the paper extends classical ergodic theorems, including Birkhoff's
theorem, to the AMS setting. It relates AMS to absolute continuity (AC)
properties for probability measures welling up from point process theory
and it generalizes results that are well known for stationary point
processes.

The general theory of AMS probability measures, with an underlying shift
transformation (say) $T$, was mainly developed during the period
1945--1985. Dowker \cite{D51} proved that, for invertible and nonsingular
$T$, the ergodic theorem holds if and only if AMS is valid. Rechard
\cite{R56} and finally Gray and Kieffer \cite{GK80} derived similar
results under
weaker assumptions for $T$. AMS is frequently used to generalize
stationarity, for instance, in information theory. See also Faigle and
Sch\"{o}nhuth \cite{FS07}.

In point process theory, it is the presence of the \textit{two} types of
shifts and \textit{two} types of stationarity which offers new
possibilities. Furthermore, in point process theory we can employ the
close relationships between the following probability distributions:
\begin{itemize}
\item[-] the distribution $P$ of the point process,

\item[-] the distributions $P_{n}$; that is, $P$ as experienced
\textit{at} the $n$th occurrence, $n\in\mathbb{Z}$,

\item[-] the distribution $P^{*}$ arising from $P_{0}$ by shifting the
origin to a completely random position between 0 and the first positive
occurrence,

\item[-] the Palm distributions $P^{x}$; that is: $P$ given an
occurrence at $x$, $x\in\mathbb{R}$,
\item[-] the shifted Palm distributions $P^{0,x}$; that is, $P$ as
experienced \textit{at} an occurrence in $x$, $x\in\mathbb{R}$.
\end{itemize}

Some researchers noted the close connection between (Birkhoff's) ergodic
results and the AMS concepts in point process theory. Daley and
Vere-Jones (\cite{DV08}; Chapter 13) give overviews of authors and
results; see
also Sigman (\cite{S95}; Chapter 2). They use coupling results to
study AMS, as in
Thorisson \cite{T00}. However, to the best knowledge of the author of the
present paper, a precise and rather complete study of AMS for point
processes on $\mathbb{R}$ has not yet been performed.

Starting point of the paper is the theorem which states that for $P$ the
concepts `event-asymptotic mean stationarity' (EAMS) and
`time-asymptotic mean stationarity' (TAMS) are equivalent. This result
and relationships between the (above mentioned) distributions are used
to link AMS to several AC properties. Thus, the class of AMS point
processes is characterized in many ways and well-known results from
stationary point process theory are generalized.

In the remainder of the (current) Section \ref{sec1}, basic notations and
definitions are introduced. In Section \ref{sec2}, we summarize
important results
of stationary point process theory. Especially the so-called
Radon--Nikodym (RN) approach -- typical for this research -- is
explained. Many of the results will be generalized later under (weaker)
AMS conditions. In Section \ref{sec3}, some well-known formulae from
\textit{nonstationary} point process theory on $\mathbb{R}$ are
rewritten into formulae resembling results from Section \ref{sec2},
into forms
useful for later sections. General definitions of $P^{x}$, $P^{0,x}$,
$P_{n}$ and $P^{*}$ are given. In Section \ref{sec4}, the concepts
TAMS and EAMS
are defined and their equivalence is proven. Also the relationship with
accompanying ergodic results is considered. Section \ref{sec5} is
about the
equivalence of AMS to AC properties for $\{P_{n}\}$ and $P^{*}$, and to
a weak AC property for $\{P^{0,x}\}$. Results from Section \ref{sec2} are
generalized by using results from Section~\ref{sec3}. Sections~\ref
{sec6} and~\ref{sec7} are about
AC properties for $\{P^{0,x}\}$ and $P$ -- respectively, with respect to
an event-stationary and a time-stationary distribution --, both stronger
than AMS. Again, results from Section \ref{sec2} are generalized. In
Section \ref{sec8},
AC properties for $P$ and $\{P^{0,x}\}$ are related and the
relationships between the RN-derivatives are considered.

\subsection*{Basic notations}

In the present research, $\mathbb{R}$ denotes the set of real
numbers and $\operatorname{Bor}(\mathbb{R}$) the set of Borel-sets of
$\mathbb{R}$. For $k\in\mathbb{Z}$, the set $\mathbb{R}_{k}$ is
defined as the positive half-line $(0,\infty)$ if $k > 0$ and
as the nonpositive half-line $(-\infty,0]$ if $k\leq0$. The notations
$:=$ and $\Leftrightarrow_{\mathrm{def}}$ both mean \textit{is by
definition}. Furthermore, a.e. means \textit{almost everywhere} and
w.r.t. means \textit{with respect to}.

Although many results in this paper can be generalized to more general
(like marked) point processes, we will only consider point processes on
$\mathbb{R}$. A \textit{point process} is a measurable mapping
$\Phi$ from a probability space $(\Omega,\mathcal{A},\mathbb{P})$
to the set $E$ of all integer-valued measures $\varphi$ on
$\mathbb{R}$ for which $\varphi(B)<\infty$ for all bounded
$B\in\operatorname{Bor}(\mathbb{R})$. $E$ is endowed with the
$\sigma
$-field $\mathcal{E}$ generated by the sets $\{\varphi\in E\dvt
\varphi(B)=k\}$, for $k\in\mathbb{Z}$ and sets
$B\in\operatorname{Bor}(\mathbb{R})$. We also define
\[
M:=\bigl\{\varphi\in E\dvt\varphi(\mathbb{R})>0; \varphi\{y\}\leq
1 \mbox{ for
all } y\in\mathbb{R}\bigr\},
\]
and add the $\sigma$-field
$\mathcal{M}:=M\cap\mathcal{E}$. We
denote the (probability) distribution of $\Phi$ by $P$. Reversely,
probability distributions on $(M,\mathcal{M})$ are
called \textit{point process distributions}. We will only allow single
occurrences; we assume that $P(M)=1$.

The atoms (called \textit{points}, \textit{events},
\textit{occurrences}, \textit{arrivals}) of $\varphi\in M$ are denoted
by $T_{n}(\varphi)$ under the convention that
\[
\cdots<T_{-1}(\varphi)<T_{0}(\varphi)\leq
0<T_{1}(\varphi)<T_{2}(\varphi)<\cdots,
\]
provided that they are finite. Occasionally, we will also write $T(n)$
instead of $T_{n}$. We write
$\alpha_{n}(\varphi):=T_{n+1}(\varphi)-T_{n}(\varphi)$,
$n\in\mathbb{Z}$, for the \textit{interval lengths} between finite
occurrences. Sets in $\mathcal{M}$ will be called
\textit{eventualities}. Eventualities like the set $\{\varphi\in
M\dvt\alpha_{n}(\varphi)<3\}$ will shortly be written as
$[\alpha_{n}(\varphi)<3]$ or even $[\alpha_{n}<3]$. Some other subsets
of $M$ with natural $\sigma$-fields:
\begin{eqnarray*}
&&F_{n}:=\bigl\{\varphi\in M\dvt\big|T_{n}(\varphi)\big|<\infty
\bigr\}\quad\mbox{and} \quad\mathcal{F}_{n}:=F_{n}\cap
\mathcal{M},\qquad n\in\mathbb{Z},
\\
&&M_{x}:=\bigl\{\varphi\in M\dvt\varphi\{x\}=1\bigr\}\quad\mbox{and}
\quad\mathcal{M}_{x}:=M_{x}\cap\mathcal{M},\qquad x\in
\mathbb{R},
\\
&&M^{\infty} :=\bigl\{\varphi\in M\dvt\varphi(-\infty,0]=\varphi
(0,\infty)=
\infty\bigr\}\quad\mbox{and} \quad\mathcal{M}^{\infty
}:=M^{\infty}
\cap\mathcal{M},
\\
&&M^{0}:=\bigl\{\varphi\in M^{\infty} \dvt\varphi\{0\}=1\bigr\}
\quad\mbox{and} \quad\mathcal{M}^{0}:=M^{0}\cap\mathcal{M}.
\end{eqnarray*}

The family $\{\theta_{y}\dvt y\in\mathbb{R}\}$ of \textit{time-shifts}
$\theta_{y} \dvtx E\rightarrow E$ defined by
$\theta_{y}(\varphi):=\theta_{y}\varphi:=\varphi(y+\cdot)$
is important. The same holds for the family
$\{\eta_{n}\dvt n\in\mathbb{Z}\}$ of \textit{event-shifts} $\eta
_{n} \dvtx F_{n}\rightarrow E$ with
$\eta_{n}(\varphi):=\eta_{n}\varphi:=\varphi(T_{n}(\varphi)+\cdot)$.
Note that $\theta_{y}\circ\theta_{x}=\theta_{y+x}$ for all $y, x\in
\mathbb{R}$ and that $\theta_{y}\varphi$ has occurrences in
$T_{k}(\varphi)-y$ (if finite) for $k\in\mathbb{Z}$. Also note that
$\eta_{n}\circ\eta_{k}=\eta_{n+k}$ for all $n, k\in\mathbb{Z}$,
that $\eta_{n}\varphi$ has occurrences in
$T_{k}(\varphi)-T_{n}(\varphi)$ (if finite), and that
$\eta_{m}=(\eta_{1})^{m}$ for all positive $m\in\mathbb{Z}$.
Regarding these shifts, the following notations are adopted:
\begin{eqnarray*}
&&\theta_{y}^{-1}A:=\{\varphi\in E\dvt\theta_{y}
\varphi\in A\},\qquad y\in\mathbb{R}\mbox{ and } A\in\mathcal{E},
\\
&&\eta_{n}^{-1}A:=\{\varphi\in F_{n}\dvt
\eta_{n}\varphi\in A\},\qquad n\in\mathbb{Z} \mbox{ and } A\in
\mathcal{E},
\\
&&\mathcal{I}':=\bigl\{A\in\mathcal{M}^{\infty} \dvt
\theta_{y}^{-1}A=A \mbox{ for all } y\in\mathbb{R}\bigr\}\quad
\mbox{and}\quad\mathcal{I}:=\bigl\{A\in\mathcal{M}^{\infty} \dvt
\eta_{1}^{-1}A=A\bigr\}.
\end{eqnarray*}
In Nieuwenhuis (\cite{N94}; Lemma 2), it was proved that the invariant
$\sigma
$-fields $\mathcal{I}'$ and $\mathcal{I}$ coincide. As a consequence,
it holds for all $\mathcal{I}$-measurable functions $f \dvtx
M^{\infty}
\rightarrow\mathbb{R}$ that:
%
\begin{equation}
\label{eq1.1} f\circ\theta_{y}=f\quad\mbox{and}\quad f\circ
\eta_{n}=f\qquad\mbox{for all }y\in\mathbb{R}\mbox{ and }n\in
\mathbb{Z}.
\end{equation}
For $A\in \mathcal{M}$ and $B\in\operatorname{Bor}(\mathbb{R})$,
we define: an \textit{A-occurrence} is an arrival time $T_{n}$ for
which the eventuality $[\eta_{n}\varphi\in A]$ occurs, $N(B)$ is the
number of occurrences in $B$ and $N_{A}(B)$ the number of the
$A$-occurrences in $B$. That is:
%
\begin{equation}
N(B):=\sum_{n \in\mathbb{Z}}1_{[T(n)\in B]}\quad\mbox{and}
\quad N_{A}(B):=\sum_{n \in\mathbb{Z} }(1_{[T(n)\in
B]}1_{A}
\circ\eta_{n}).
\end{equation}
Expectation under $\mathbb{P}$ is denoted by
$\mathbb{E}$, expectation under $P$ by $E$. For
measurable functions $f\dvtx M\rightarrow\mathbb{R}$, we use
$\mathbb{E}f(\Phi)$, $Ef$, $Ef(\Phi)$ and even $Ef(\varphi)$ to
denote the expectation of~$f(\Phi)$.

From Section \ref{sec4} onwards, we will use the notations $P_{ts}$
and $P_{es}$
to respectively denote an event-stationary (ES) and a time-stationary
(TS) distribution on $(M,\mathcal{M})$. Furthermore, AC means `absolute
continuity'. For two probability measures $P$ and $Q$ on
$(M,\mathcal{M})$, the notation $P\ll Q$ denotes that $P$ is AC with
respect to (w.r.t.) $Q$. We also say that $Q$ \textit{dominates} $P$
and denote a
Radon--Nikodym derivative as RN. A discrete-time stochastic process
$\{X_{n}\dvt n\in\mathbb{Z}\}$ on $M$ is called $Q$-\textit{stationary}
or \textit{stationary under} $Q$ if it holds for each positive integer $m$
that:
%
\begin{equation}
(X_{1}, \dots, X_{m})=_{d}(X_{k+1},
\dots, X_{k+m})\qquad\mbox{under }Q\ (\mbox{for all }k\in\mathbb{Z}).
\end{equation}

\section{Stationary point processes}\label{sec2}

We offer a brief but coherent overview of stationary point process
theory, enclosing only results that will be used or generalized later;
our notations originate from Franken \textit{et al.} \cite{FKAS82}.
The second half of
this section is less known; it reflects the special approach in the
present paper.

A point process $\Phi$ (and also its distribution $P$) is called
\textit{time-stationary} (shortly TS) if
$P\theta_{y}^{-1}(A):=P(\theta_{y}^{-1}A)=P(A)$ for all
$y\in\mathbb{R}$ and $A\in \mathcal{M}$. It is called
\textit{event-stationary} (ES) if $P(M^{\infty} )=1$ and it holds for
all $A\in\mathcal{M}^{\infty}$ that
$P\eta_{1}^{-1}(A):=P(\eta_{1}^{-1}A)=P(A)$, and hence that
$P(\eta_{n}^{-1}A)=P(A)$ for all $n\in\mathbb{Z}$. For TS
distributions $P$ with $P(M)=1$, we have $P(M^{\infty} )=1$ and
$P(M^{0})=0$; ES distributions $P$ satisfy$ P(M^{\infty} )=1$ and
$P(M^{0})=1$. If $\Phi$ is TS, we call $\lambda:=E(N(0,1])$ the
\textit{intensity} of $\Phi$ and $P$; we will always implicitly assume
that \mbox{$\lambda<\infty$}.

Suppose that $P$ is TS and $y \geq0$. Then, $E(N(0,y])=\lambda y$.
For all $x > 0$ the definition below yields one probability measure
$P^{0}$ on $(M^{\infty},\mathcal{M}^{\infty} )$, the \textit{Palm
distribution} (PD) of $\Phi$ and $P$:
%
\begin{equation}
\label{eq2.1} P^{0}(A):=\frac{1}{\lambda x}E\bigl(N_{A}(0,x]\bigr)=
\frac{1}{\lambda
x}E\Biggl(\sum_{i=1}^{N(0,x]}1_{A}
\circ\eta_{i}\Biggr)\qquad\mbox{for } A\in\mathcal{M}^{\infty}.
\end{equation}
Informally, $P^{0}$ is the conditional distribution of the point process
if there is an occurrence in the origin. We denote $P^{0}$-expectation
by $E^{0}$. This PD has the following properties:
%
\begin{eqnarray}
&&P^{0}\bigl(M^{0}\bigr)=1,\qquad P^{0}
\eta_{n}^{-1}=P^{0}\qquad\mbox{for all }n\in
\mathbb{Z},
\\
&&\lambda=\frac{1}{E^{0}(\alpha_{0})}=E\biggl(\frac{1}{\alpha
_{0}}\biggr),
\\
\label{eq2.4} &&P^{0} (A )=\frac{E(N_{A}(0,x])}{E(N(0,x])}\qquad
\mbox{for all } x>0
\mbox{ and }A\in\mathcal{M}^{\infty}.
\end{eqnarray}
Hence, the PD of a TS distribution $P$ is ES and the sequence
$\{\alpha_{n}\}$ is stationary under $P^{0}$. With
$\lambda_{A}:=E(N_{A}(0,1])$, the \textit{intensity of the
A-occurrences}, it follows that:
%
%
\begin{equation}
\label{eq2.5} P^{0} (A )=\lambda_{A}/\lambda.
\end{equation}
Compared to (\ref{eq2.1}), the following so-called \textit{inversion formulae}
work the other way round:
%
\begin{eqnarray}
\label{eq2.6} P(A) &=&\lambda\int_{\mathbb{R}_{k}}P^{0}\bigl
[\varphi(-x+
\cdot)\in A \mbox{ and } T_{-k}(\varphi)\leq-x<T_{-k+1}(
\varphi)\bigr]\,\mathrm{d}x
\nonumber
\\[-8pt]
\\[-8pt]
&=&\lambda E^{0} \biggl(\int_{-T(-k+1)}^{-T(-k)}1_{A}
\circ\theta_{-x}\,\mathrm{d}x \biggr)=\lambda E^{0} \biggl(
\int_{T(-k)}^{T(-k+1)}1_{A}\circ
\theta_{y}\,\mathrm{d}y \biggr).
\nonumber
\end{eqnarray}
Here, $A\in\mathcal{M}^{\infty}$ and $k\in\mathbb{Z}$. It is
allowed to replace $\mathbb{R}_{k}$ by $\mathbb{R}$. See
Slivnyak \cite{S62} and Kaplan \cite{K55} for the one-to-one correspondence
described in (\ref{eq2.4}) and (\ref{eq2.6}).

In Nieuwenhuis (\cite{N89}; Theorem 8.1) it was proved that, for TS
distributions
$P$ and all $n\in\mathbb{Z}$, the \textit{intermediate
distribution} $P_{n}:=P\eta_{n}^{-1}$ is equivalent (i.e., AC
in two directions) to the PD $P^{0}$:\vspace*{-12pt}
%
\begin{subequation}\label{eq2.7}
\begin{eqnarray}
\label{eq2.7a} P_{n}&\ll& P^{0}\quad\mbox{and}\quad
P_{n}(A)=\lambda E^{0}(\alpha_{-n}1_{A}),
\\
\label{eq2.7b} P^{0}&\ll& P_{n}\quad\mbox{and}\quad
P^{0}(A)=\frac{1}{\lambda} E_{n}\biggl(\frac{1}{\alpha_{-n}}1_{A}
\biggr)=\frac{1}{\lambda} E\biggl(\frac{1}{\alpha_{0}}1_{A}\circ
\eta_{n}\biggr),\qquad A\in\mathcal{M}^{\infty}.\hspace*{8pt}
\end{eqnarray}
\end{subequation}
(We write $E_{n}$ for $P_{n}$-expectation.) See also Ryll--Nardzewski
\cite{R61} and Thorisson \cite{T00} for similar approaches. Results
(\ref{eq2.7a}), (\ref{eq2.7b}),
which reflect the so-called \textit{Radon--Nikodym} \textit{approach},
offer the opportunity to jump easily between $P$, $P^{0}$ and related
distributions and are very important for this paper. We illustrate their
use and derive some frequently used results.

Since $ P(A) =\lambda E^{0}(\int_{0}^{\alpha_{0}}1_{A}\circ\theta
_{y}\,\mathrm{d}y)$,
it follows from (\ref{eq2.7b}) that $P(A)$ can be written otherwise as a
$P$-expectation:
%
\begin{subequation}\label{eq2.8}
\begin{equation}
\label{eq2.8a} P(A)=E_{0}\biggl(\frac{1}{\alpha_{0}}\int
_{0}^{\alpha
_{0}}1_{A}
\circ\theta_{y}\,\mathrm{d}y\biggr) = E\biggl(\frac{1}{\alpha
_{0}}\int
_{T_{0}}^{T_{1}}1_{A}\circ\theta_{y}
\,\mathrm{d}y\biggr).
\end{equation}
By (\ref{eq2.8a}) we obtain, for all functions $g\dvtx\mathbb
{R}\rightarrow\mathbb{R}$ with $E|g(T_{1})|<\infty$:
\begin{equation}
\label{eq2.8b} E\bigl(g(T_{1}) | (\alpha_{n})_{n \in\mathbb{Z}}
\bigr)=\frac{1}{\alpha_{0}}\int_{0}^{\alpha_{0}}g(x)\,
\mathrm{d}x \qquad P\mbox{-a.s.}
\end{equation}
Hence: conditionally on $\alpha_{0}$, the distribution of $T_{1}$ under
$P$ is $\operatorname{uniform}(0, \alpha_{0})$. Note that
$h:=\frac{1}{\alpha_{0}}\int_{T_{0}}^{T_{1}}1_{A}\circ\theta_{y}\,
\mathrm{d}y$ satisfies
$h\circ\eta_{0}=h$ on $M^{\infty}$. We obtain by (\ref{eq2.8a}), for
$P$-integrable functions $f$ and $g\dvtx M^{\infty} \rightarrow
\mathbb{R}$:
\begin{equation}
\label{eq2.8c} E\biggl(f\cdot\frac{1}{\alpha_{0}}\int_{T(0)}^{T(1)}g
\circ\theta_{y}\,\mathrm{d}y\biggr) = E\biggl(g\cdot\frac
{1}{\alpha_{0}}
\int_{T(0)}^{T(1)}f\circ\theta_{y}\,
\mathrm{d}y\biggr).
\end{equation}
\end{subequation}
Set $\overline{\alpha} =E^{0}(\alpha_{0}\vert\mathcal{I})$ and
$\overline{N}=E(N(0,1]|\mathcal{I})$. By (\ref{eq2.7a}) and (\ref{eq2.7b}),
Birkhoff's ergodic results
\[
\frac{1}{n}\sum_{i=1}^{n}
\alpha_{i}\rightarrow\overline{\alpha}\qquad\mbox{as }
n\rightarrow
\infty~P^{0}\mbox{-a.s.}\quad\mbox{and}\quad\frac{1}{x}N(0,x]
\rightarrow\overline{N}\qquad\mbox{as }x\rightarrow\infty~ P\mbox{-a.s.}
\]
are also valid $P$-a.s. and $P^{0}$-a.s., respectively. Furthermore, it
can be proved that:
%
\begin{equation}
\overline{N}=\frac{1}{\overline{\alpha}} =E\biggl(\frac{1}{\alpha
_{0}}\Big\vert\mathcal{I}\biggr)
\qquad P^{0}\mbox{-a.s. and }P\mbox{-a.s.}
\end{equation}
A TS point process is called \textit{ergodic} if $P(C)=0$ or 1 (or,
equivalently, $P^{0}(C)=0$ or 1) for all
$C\in\mathcal{I}$. It is called \textit{pseudo-ergodic} if
$P^{0}[\lambda\overline{\alpha} =1]=1$; see also Nieuwenhuis \cite{N98}.

Note that, for all $x\in\mathbb{R}$, $A\in\mathcal{M}^{\infty}$
and functions $f\dvtx M^{\infty} \rightarrow\mathbb{R}$ with
$f=f\circ\eta_{0}$ on $M^{\infty}$:
\begin{subequation}\label{eq2.10}
\begin{eqnarray}
\label{eq2.10a} E^{0}(1_{A}\cdot f\circ\theta_{-x})&
= &\sum_{k\in\mathbb{
Z}}E^{0}(1_{A}\cdot f
\circ\eta_{k}\cdot1_{[T_{k} \leq -x < T_{k+1}]})
\nonumber
\\
&= &E^{0}\bigl(f\cdot N_{A}[x,x+\alpha_{0})\bigr)
\\
\label{eq2.10b} &= &\frac{1}{\lambda} E\bigl(f\cdot N_{A}[x+T_{0},x+T_{1}
)/\alpha_{0}\bigr),
\end{eqnarray}
\vspace*{-14pt}
\begin{equation}
\label{eq2.10c} P^{0}(A)=E\biggl(N_{A}[x+T_{0},x+T_{1}
)\frac{1}{\lambda\alpha_{0}}\biggr)\quad\mbox{and}\quad E\biggl
(N[x+T_{0},x+T_{1}
)\frac{1}{\alpha_{0}}\biggr)=\lambda.
\end{equation}
\end{subequation}
In coming sections, ES distributions and TS distributions are usually
denoted as $P_{es}$ and $P_{ts}$ (and the accompanying expectation
operators as $E_{es}$ and $E_{ts}$). The ES Palm distribution associated
with a TS distribution $P_{ts}$, is denoted by $P_{ts}^{0}$. So, the
relationships between $P_{ts}$ and $P_{ts}^{0}$ are the same as the
relationships between $P$ and $P^{0}$ described in (\ref
{eq2.1})--(\ref{eq2.10}).

\section{Non-stationary point processes}\label{sec3}

In this section, we consider, for \textit{general} point processes on
$\mathbb{R}$, the PDs $\{P^{x}\}$ and their shifted versions
$\{P^{0,x}\}$. Furthermore, we carefully define the distributions
$P_{n}$ and $P^{*}$ informally mentioned in Section \ref{sec1}. Campbell's
equation is used to obtain inversion formulae (\ref{eq3.6}) and (\ref
{eq3.7}) which
resemble and generalize (\ref{eq2.6}) and (\ref{eq2.7a}). We
generalize (\ref{eq2.5}) for the
case that $P=P^{*}$ and characterize the class of the TS point process
distributions.

We assume that the point process $\Phi$ satisfies $P(M)=1$, and that the
\textit{intensity measure} $\nu$ on $\operatorname{Bor}(\mathbb
{R})$ with
$\nu(B):=E(N(B))$ for $B\in\operatorname{Bor}(\mathbb{R})$ exists
and is
locally finite. Below, for $A\in \mathcal{M}$, also the
locally finite measures $\nu_{A}$ and $\mu_{A}$ play important roles:
\[
\nu_{A}(B):= E\bigl(N(B)1_{A}\bigr)\quad\mbox{and} \quad
\mu_{A}(B):=E\bigl(N_{A}(B)\bigr);\qquad B\in
\operatorname{Bor}(\mathbb{R}).
\]

\subsection*{Palm distributions}

For $A\in \mathcal{M}$, $\nu_{A}$ is dominated by $\nu$. An
RN-derivative is denoted by $x\rightarrow P^{x}(A)$, so:
%
%
\begin{equation}
\label{eq3.1} \nu_{A} (B )=\int_{B}P^{x}
(A )\,\mathrm{d}\nu(x ); \qquad B\in\operatorname{Bor} (\mathbb{R} ).
\end{equation}
A basic result in Palm theory now is that $\{P^{x}(A)\dvt x\in\mathbb
{R} \mbox{ and } A\in \mathcal{M}\}$ can
be chosen such that the function $x\rightarrow P^{x}(A)$ is measurable
for all $A\in \mathcal{M}$, $P^{x}$ is a probability measure
on $\mathcal{M}$ for all $x\in\mathbb{R}$, and
%
\begin{equation}
\label{eq3.2} \int_{M}^{}\int
_{-\infty}^{\infty} f (x,\varphi)\,\mathrm{d}\varphi(x )\,
\mathrm{d}P (\varphi)=\int_{-\infty}^{\infty} \int
_{M}^{}f (x,\varphi)\,\mathrm{d}P^{x}
(\varphi)\,\mathrm{d}\nu(x )
\end{equation}
for all $\operatorname{Bor}(\mathbb{R})\times\mathcal{M}$-measurable
functions $f$ on $\mathbb{R}\times M$ that are either nonnegative
or satisfy $\mathbb{E}[\int_{-\infty}^{\infty} f(x, \Phi) \,
\mathrm{d}\Phi(x)]<\infty$. Thus, the family $\{P^{x}\}$ of probability
distributions turns out to be uniquely defined by (\ref{eq3.2}) apart
from a
Borel-set in $\mathbb{R}$ with $\nu$-measure 0. See Matthes
\cite{M63}. See also Jagers \cite{J73} and Kallenberg \cite{K83/86}.
In the sequel,
we will assume that the family $\{P^{x}\}$ is chosen this way. Note that
$f(x,\varphi)=1_{B\times A}(x,\varphi)$, $x\in\mathbb{R}$ and
$\varphi\in M$, returns (\ref{eq3.1}). The probability measures in
$\{P^{x}\dvt x\in\mathbb{R}\}$ are called \textit{Palm distributions}
(PDs) of $P$. It can be proved that $P^{x}(M_{x})=1$ for $\nu$-a.e.
$x\in\mathbb{R}$. By letting $A$ in (\ref{eq3.1}) shrink to $\{x\}$, we
obtain the intuitive meaning for $P^{x}(A)$ as the probability that
$\Phi\in A$ under the condition that $\Phi\{x\}=1$.

\subsection*{Shifted Palm distributions}

We are especially interested in $\{P^{0,x}\}$, the family of shifted PDs
defined by $P^{0,x}:=P^{x}\theta_{x}^{-1}$. Note that
$P^{0,x}$ satisfies $P^{0,x}(M_{0})=1$, and that, in queuing terms, it
can be considered as the probability measure that under $P$ is
experienced by a customer arriving at time $x$. For time-stationary $P$
we have $P^{0,x}=P^{0}$ for $\nu$-a.e. $x\in\mathbb{R}$, where
$P^{0}$ is the event-stationary PD of $P$ in (\ref{eq2.1}). Note that
the choice
$f(x,\varphi)=1_{A}(\theta_{x}\varphi)1_{B}(x)$ in (\ref{eq3.2})
yields that:
%
\begin{equation}
\label{eq3.3} \mu_{A} (B )=\int_{B}P^{0,x}
(A )\,\mathrm{d}\nu(x ); \qquad B\in\operatorname{Bor} (\mathbb{R}
) \mbox{ and
} A\in\mathcal{M}.
\end{equation}
Hence, for all $A\in \mathcal{M}$, the function $x\rightarrow
P^{0,x}(A)$ is just an RN-derivative of $\mu_{A}$ with respect to $\nu$.
If $\nu$ is AC with respect to Leb with \textit{intensity}
$\lambda(\cdot)$, then $\mu_{A}$ is also AC with respect to Leb.
If $x\rightarrow\lambda_{A}(x)$ denotes an accompanying RN-derivative
(the \textit{intensity} of the point process of the $A$-occurrences), it
then follows for all $A\in \mathcal{M}$ that
%
\begin{equation}
\label{eq3.4} P^{0,x}(A)=\lambda_{A}(x)/\lambda(x)\qquad
\mbox{for }\nu\mbox{-a.e. } x\in\mathbb{R};
\end{equation}
cf. (\ref{eq2.5}). However, it cannot be concluded that for $\nu
$-a.e. $x\in\mathbb{R}$ the shifted PDs
satisfy $P^{0,x}(A)=\lambda_{A}(x)/\lambda(x)$ for all $A\in
\mathcal{M}$; (\ref{eq3.4}) not even necessarily defines a
probability measure
for $\nu$-a.e. $x\in\mathbb{R}$. By letting $B$
in (\ref{eq3.3}) shrink to $\{x\}$ we obtain the intuitive meaning for
$P^{0,x}(A)$ as the probability that $\theta_{x}\Phi\in A$ under the
condition that $\Phi\{x\}=1$.

\subsection*{Intermediate probability measures}

For $n\in\mathbb{Z}$ with $P(F_{n})>0$, we define the
\textit{intermediate probability measure} $P_{n}$ of $P$ as a
conditional probability distribution:
%
\begin{equation}
\label{eq3.5} P_{n}(A):=P\bigl([\eta_{n}\varphi\in A] |
F_{n}\bigr), \qquad A\in\mathcal{M}.
\end{equation}
We investigate the relationships between $P$, $\{P_{n}\}$ and
$\{P^{0,x}\}$. Set $I_{x}:=(0,x]$ for $x>0$ and
$I_{x}:=(x,0]$ for $x\leq0$, and choose $f$ in (\ref{eq3.2}) as:
\[
f(x,\varphi)=1_{A}(\varphi)1_{\{\vert
k|\}}\bigl(
\varphi(I_{x})\bigr)1_{\mathbb{R}_{k}}(x),\qquad A\in\mathcal
{M}\mbox{
and }k\in\mathbb{Z}.
\]
We obtain:
%
\begin{subequation}\label{eq3.6}
\begin{eqnarray}
\label{eq3.6a} P(A\cap F_{k})&=&\int_{\mathbb{R}_{k}}P^{x}
\bigl(A\cap\bigl[\varphi(I_{x})=|k|\bigr]\bigr)\,\mathrm{d}\nu(x),
\\
\label{eq3.6b} P(A\cap F_{k})&=&\int_{\mathbb{R}_{k}}P^{x}
\bigl(A\cap[T_{k}=x]\bigr)\,\mathrm{d}\nu(x),
\\
\label{eq3.6c} P(A\cap F_{k}) &=&\int_{\mathbb{R}_{k}}P^{0,x}
\bigl([\theta_{-x}\varphi\in A]\cap[T_{-k}\leq
-x<T_{-k+1}]\bigr)\,\mathrm{d}\nu(x).
\end{eqnarray}
\end{subequation}
Compare with (\ref{eq2.6}). It follows that, for all $A\in
\mathcal{M}$ and $k\in\mathbb{Z}$ with $P(F_{k})>0$:
%
%
\begin{equation}
\label{eq3.7} P_{k} (A )=\frac{1}{P (F_{k} )}\int_{\mathbb{R}_{k}}P^{0,x}
\bigl(A\cap[T_{-k}\leq-x<T_{-k+1} ] \bigr)\,\mathrm{d}\nu(x
).
\end{equation}
Note that it is allowed to replace $\mathbb{R}_{k}$ by
$\mathbb{R}$ in (\ref{eq3.6b}), (\ref{eq3.6c}) and (\ref{eq3.7}),
and that (\ref{eq3.7})
generalizes (\ref{eq2.7a}). Substitution of $A\cap[T_{k}\in B]$ for
$A$ in
(\ref{eq3.6c}) yields
%
\begin{equation}
\label{eq3.8} P \bigl(A\cap[T_{k}\in B ] \bigr)=\int
_{B}P^{0,x}\bigl( [\theta_{-x}\varphi\in A
]\cap[T_{-k}\leq-x<T_{-k+1} ]\bigr)\,\mathrm{d}\nu(x )
\end{equation}
for all $k\in\mathbb{Z}$, $B\in\operatorname{Bor}(\mathbb{R})$ and
$A\in \mathcal{M}$. By taking $\sum_{k\in\mathbb{Z}}$, the
left-hand side becomes equal to $\nu_{A}(B)$ and we get (\ref{eq3.1})
back. When $A$ in (\ref{eq3.8}) is replaced by $[\eta_{k}\varphi\in
A]$, we
obtain:
%
%
\begin{equation}
\label{eq3.9} P \bigl([\eta_{k}\varphi\in A]\cap[T_{k}
\varphi\in B ] \bigr)=\int_{B}P^{0,x}\bigl(A\cap
[T_{-k}\leq-x<T_{-k+1} ]\bigr)\,\mathrm{d}\nu(x ).
\end{equation}
Note that we get (\ref{eq3.3}) back by taking $\sum_{k\in\mathbb{Z}}$.
The choice $A=M_{0}$ in (\ref{eq3.9}) ensures that, if $P(F_{k})$ is
larger than
0, the conditional distribution $P([T_{k}\in\cdot] \vert F_{k})$ of
$T_{k}$ is AC with respect to $\nu$ with RN-derivative
$\gamma(x):=P^{0,x}[T_{-k}\leq-x<T_{-k+1}]/P(F_{k})$. So:
%
\begin{equation}
\label{eq3.10} P\bigl([T_{k}\in\cdot] | F_{k}\bigr)\ll\nu,
\qquad\gamma(x)=P^{x}[T_{k}=x]/P(F_{k})\qquad
\mbox{for }\nu\mbox{-a.e. }x\in\mathbb{R}.
\end{equation}

\subsection*{The distribution $P^{*}$}

For $P$ such that $P(M^{\infty} )=1$, we define the distribution $P^{*}$
as follows:
%
\begin{equation}
\label{eq3.11} P^{*}(A):=E\biggl(\frac{1}{\alpha_{0}}\int
_{T(0)}^{T(1)}1_{A}\circ\theta_{y}
\,\mathrm{d}y\biggr)=E_{0}\biggl(\frac{1}{\alpha_{0}}\int
_{0}^{\alpha_{0}}1_{A}\circ\theta_{y}
\,\mathrm{d}y\biggr)\qquad\mbox{for } A\in\mathcal{M}^{\infty}.
\end{equation}
By (\ref{eq2.8a}) and (\ref{eq2.7a}), time-stationary point processes
(with finite
intensity) satisfy:
\begin{enumerate}[(b)]
\item[(a)] $P=P^{*}$\quad and
%
\begin{subequation}\label{eq3.12}\vspace*{-20pt}
\begin{equation}
\label{eq3.12a}
\end{equation}
\item[(b)] there exists an ES point process distribution $P_{es}$ such
that:
\begin{equation}
\label{eq3.12b} P_{0}\ll P_{es}\quad\mbox{and}\quad
\mathrm{d}P_{0}/\mathrm{d}P_{es}=\lambda\alpha_{0}\qquad\mbox
{with }
\lambda=1/E_{es}(\alpha_{0}) \in(0,\infty).
\end{equation}
\end{subequation}
\end{enumerate}
Reversely, if $P$ satisfies (\ref{eq3.12a}), (\ref{eq3.12b}) then, for
$A\in\mathcal{M}^{\infty}$:
\[
P(A)=P^{*}(A)=\lambda E_{es}\biggl(\int_{0}^{\alpha
_{0}}1_{A}
\circ\theta_{y}\,\mathrm{d}y\biggr)=:P_{ts}(A).
\]
Note that $P_{ts}$ (and hence $P$) is just the TS distribution such that
the accompanying PD is $P_{es}$; see (\ref{eq2.6}). Hence, (\ref
{eq3.12a}), (\ref{eq3.12b})
characterize the class of TS point process distributions.

For $A\in\mathcal{M}^{\infty}$ and $B\in\operatorname{Bor}(\mathbb{R})$,
set $\mu_{A}^{*}(B):=E^{*}(N_{A}(B))$ and $\nu^{*}(B):=E^{*}(N(B))$.
Here, $ E^{*}$~refers to $P^{*}$-expectation.

\begin{theorem}\label{teo3.1}
Suppose that $P$ is a point process distribution
with $P(M^{\infty} )=1$. Then:
\begin{enumerate}[(4)]
\item[(1)] the intermediate distributions of $P^{*}$ and $P$ coincide;

\item[(2)] under $P^{*}$, the conditional distribution of $T_{1}$ given
$(\alpha_{n})_{n \in\mathbb{Z}}$, is $\operatorname{uniform}(0,
\alpha_{0})$;

\item[(3)] for $A\in\mathcal{M}^{\infty}$ it holds that
$\mu_{A}^{*}\ll\mathrm{Leb}$ and $\nu^{*}\ll\mathrm{Leb}$, with
intensity functions
\[
\lambda_{A}^{*}(x)=E\biggl(\frac{1}{\alpha_{0}}N_{A}[x+T_{0},x+T_{1}
)\biggr)\quad\mbox{and}\quad\lambda^{*}(x)=E\biggl(
\frac{1}{\alpha_{0}}N[x+T_{0},x+T_{1})\biggr);
\]
\item[(4)] the shifted PDs of $P^{*}$ satisfy
$P^{*0,x}(A)=\lambda_{A}^{*}(x)/\lambda^{*}(x)$ for $\nu^{*}$-a.e.
$x\in\mathbb{R}$ and \mbox{$A\in\mathcal{M}^{\infty}$}.
\end{enumerate}
\end{theorem}

\begin{pf}
Part (1) is immediate. Part (2) holds since for all
eventualities $A$ in the $\sigma$-field generated by
$(\alpha_{n})_{n \in\mathbb{Z}}$ and all functions $g\dvtx\mathbb
{R}\rightarrow R$ with $E^{*}|g(T_{1})|<\infty$ we have:
\begin{eqnarray*}
E^{*}\bigl(1_{A}E^{*}\bigl(g(T_{1})|(
\alpha_{n})_{n \in\mathbb{Z}}\bigr)\bigr)& =& E^{*}
\bigl(1_{A}g(T_{1})\bigr) = E\biggl(1_{A}
\frac{1}{\alpha_{0}}\int_{T(0)}^{T(1)}g(T_{1}
\circ\theta_{y})\,\mathrm{d}y\biggr)
\\
&= &E_{0}\biggl(1_{A}\frac{1}{\alpha_{0}}\int
_{0}^{\alpha_{0}}g(x)\,\mathrm{d}x\biggr) = E^{*}
\biggl(1_{A}\frac{1}{\alpha_{0}}\int_{0}^{\alpha_{0}}g(x)
\,\mathrm{d}x\biggr).
\end{eqnarray*}
For (3), note that for $A\in\mathcal{M}^{\infty}$ and
$B\in\operatorname{Bor}(\mathbb{R})$ we have that $\mu_{A}^{*}(B)$
equals:
\begin{eqnarray*}
&&\sum_{k\in\mathbb{Z}}P^{*}\bigl[T_{k}(
\varphi)\in B \mbox{ and } \eta_{k}\varphi\in A\bigr]
\\
&&\quad= \sum_{k\in\mathbb{Z}}E\biggl(\frac{1}{\alpha_{0}}\int
_{T(k)-T(1)}^{T(k)-T(0)}1_{B}(y)\,\mathrm{d}y1_{[\eta_{k}\varphi
\in A]}
\biggr)
\\
&&\quad= \sum_{k\in\mathbb{Z}}E\biggl(\frac{1}{\alpha_{0}}\int
_{B}1_{[y+T(0) \leq
T(k) < y+T(1)]}\,\mathrm{d}y1_{[\eta_{k}\varphi \in A]}\biggr) =
\int_{B}\lambda_{A}^{*}(x)\,\mathrm{d}x.
\end{eqnarray*}
Hence, $\mu_{A}^{*}\ll\mathrm{Leb}$. The choice $A=M^{\infty}$ yields
$\nu^{*}\ll\mathrm{Leb}$. For (4), first note that it holds for
$\nu^{*}$-a.e. $x\in\mathbb{R}$ that $Q^{0,x}(A):= \lambda
_{A}^{*}(x)/\lambda^{*}(x)$ defines a probability measure on
$(M^{\infty},\mathcal{M}^{\infty} )$. Replacing $P$, $\{P^{x}\}$ and
$\nu$ by $P^{*}$, $\{Q^{0,x}\theta_{-x}^{-1}\}$ and $\nu^{*}$ in
(\ref{eq3.2})
yields, for both sides:
\[
E\biggl(\frac{1}{\alpha_{0}(\varphi)}\sum_{k \in
\mathbb{Z}}\int
_{T_{0}\varphi}^{T_{1}\varphi} f(T_{k}\varphi-y,
\theta_{y}\varphi)\,\mathrm{d}y\biggr).
\]
So, (4) follows from the uniqueness of the family of PDs.
\end{pf}

It follows from Theorem \ref{teo3.1}(1), (2) that $P^{*}$ arises from
$P_{0}$ by
shifting the origin to an arbitrary position in the interval
$(0,\alpha_{0})$. By (\ref{eq2.10c}) and Theorem \ref{teo3.1}(4), we
get a generalization of
(\ref{eq2.5}); cf. (\ref{eq2.10c}):
%
\begin{equation}
\label{eq3.13} P=P^{*}\quad\Rightarrow\quad P^{0,x}(A)=
\lambda_{A}(x)/\lambda(x)\qquad\mbox{for } \nu\mbox{-a.e. } x\in
\mathbb{R}\mbox{ and all }A\in\mathcal{M}^{\infty}.
\end{equation}

\section{Asymptotic mean stationarity}\label{sec4}

A point process (as well as its distribution $P$) is called
\textit{time-asymptotic(ally) mean stationary} (shortly TAMS) if a
probability distribution $P_{ts}$ on $(M,\mathcal{M})$ exists such that:
%
%
\begin{equation}
\label{eq4.1} \frac{1}{x}\int_{0}^{x}P [
\theta_{y}\varphi\in A ]\,\mathrm{d}y \rightarrow P_{ts} (A )
\qquad\mbox{as } x\rightarrow\infty, \mbox{ for all } A\in
\mathcal{M}.
\end{equation}
Note that $P_{ts}$ is indeed TS. We write ``$P$ is TAMS($P_{ts}$)'' and
call $P_{ts}$ the \textit{time-stationary} \textit{limit distribution}
of $P$. A point process (and its distribution $P$) with $P(M^{\infty}
)=1$ is called \textit{event-asymptotic(ally) mean stationary} (shortly
EAMS) if a probability distribution $P_{es}$ on $(M^{\infty}
,\mathcal{M}^{\infty} )$ exists such that:
%
%
\begin{equation}
\label{eq4.2} \frac{1}{n}\sum_{i=1}^{n}P
[\eta_{i}\varphi\in A ]\rightarrow P_{es} (A )\qquad
\mbox{as } n\rightarrow\infty, \mbox{ for all } A\in\mathcal
{M}^{\infty}.
\end{equation}
We write ``$P$ is EAMS($P_{es}$)'' and call $P_{es}$ the
\textit{event-stationary} \textit{limit distribution} of $P$. Note that
$P_{es}$ is ES and that, for $P$ being EAMS, it is only its behavior on
$\mathcal{M}^{0}$ which matters. Sigman \cite{S95} refers to the AMS
concepts as time and event asymptotic stationarity; Daley and Vere-Jones
\cite{DV08} use $(C, 1)$-asymptotic stationarity and event-stationarity,
respectively.

Note that TS point processes are TAMS and ES point processes are EAMS.
However, we will see that the class of TAMS (EAMS) point processes is --
considerably -- larger than the class of TS (ES) point processes. As an
example: with $P_{ts}$ the distribution of the TS Poisson point process
with intensity $\lambda_{ts}$, the definition
$P(A):=E_{ts}(1_{A}\cdot N(0, 1])/\lambda_{ts}$ for
$A\in\mathcal{M}^{\infty}$ yields a point process distribution $P$ which
is absolutely continuous with respect to $P_{ts}$ and hence (see (\ref{eq4.8b})
below) is TAMS. However, $P$ is \textit{not} TS since, by Jensen's
inequality:
\[
E\bigl(N(0, 1]\bigr)=E_{ts}\bigl[\bigl(N(0, 1]\bigr)^{2}\bigr]/
\lambda_{ts}>\lambda_{ts}= E\bigl(N(1, 2]\bigr).
\]

The following characterizations of EAMS and TAMS will be used
frequently.

\begin{theorem}\label{teo4.1}
Let $\Phi$ be a point process with distribution $P$
for which $P(M^{\infty} )=1$, and let $P_{es}$ and $ P_{ts}$,
respectively, be an ES and a TS point process distribution. Then:
%
\begin{eqnarray}
\label{eq4.3} P\mbox{ is }\operatorname{EAMS}(P_{es})\quad&
\Leftrightarrow&\quad P=P_{es}\qquad\mbox{on }\mathcal{I};
\nonumber
\\[-8pt]
\\[-8pt]
P \mbox{ is }\operatorname{TAMS}( P_{ts}) \quad&\Leftrightarrow
&\quad
P=P_{ts}\qquad\mbox{on }\mathcal{I}.
\nonumber
\end{eqnarray}
\end{theorem}

\begin{pf}
The implications `$\Rightarrow$' follow from (1.1). For
the implications `$\Leftarrow$', suppose respectively that $P=P_{es}$
on $\mathcal{I}$ and $P=P_{ts}$ on $\mathcal{I}$. Note that, for each
$A\in\mathcal{M}^{\infty}$, the eventualities
\[
\Biggl[\frac{1}{n}\sum_{i=1}^{n}1_{A}
\circ\eta_{i}\rightarrow E_{es}(1_{A}|\mathcal{I})
\Biggr]\quad\mbox{and}\quad\biggl[\frac{1}{x}\int_{0}^{x}1_{A}
\circ\theta_{y}\,\mathrm{d}y\rightarrow E_{ts}(1_{A}|
\mathcal{I})\biggr]
\]
are elements of $\mathcal{I}$ which (by Birkhoff's ergodic theorem) have
probability 1 under respectively, $P_{es}$ and $P_{ts}$, and hence they
both have $P$-probability 1. Take $P$-expectations.
\end{pf}

The next theorem roughly states that $P$ is EAMS iff $P$ is TAMS; see
also Daley and Vere-Jones (\cite{DV08}; Theorem 13.4.VI) for a similar (but
different) theorem. If $P$ is $\operatorname{EAMS}(P_{es})$, then
$\lim_{n\rightarrow
\infty} \sum_{i=1}^{n}\alpha_{i}/n$ equals $\overline{\alpha}
:=E_{es}(\alpha_{0}\vert\mathcal{I})$; this holds $P_{es}$-a.s. and (by
Theorem \ref{teo4.1}) also \mbox{$P$-a.s.} If $P$ is $\operatorname
{TAMS}(P_{ts})$, then
$\lim_{x\rightarrow\infty} N(0,x]/x$ equals
$\overline{N}:= E_{ts}(N(0,1]|\mathcal{I})$, $ P_{ts}$-a.s. and also
\mbox{$P$-a.s.}

\begin{theorem}\label{teo4.2}
Let $\Phi$ be a point process with distribution $P$
for which $P(M^{\infty} )=1$. Then:
\begin{eqnarray*}
&&P\mbox{ is }\operatorname{EAMS}(P_{es})\quad\mbox{and} \quad
P_{es}[0<\overline{\alpha} <\infty]=1
\\
&&\hspace*{86pt}\Leftrightarrow\\
&&P\mbox{ is }\operatorname{TAMS}(P_{ts})
\quad\mbox{and}\quad P_{ts}[0<\overline{N}<\infty]=1.
\end{eqnarray*}
\noindent These distributions $P_{es}$ and $P_{ts}$ are related as follows:
%
\begin{eqnarray}
\label{eq4.4} &&P_{ts}(A)=E_{es}\biggl(\frac{1}{\overline{\alpha}}
\int_{0}^{\alpha
_{0}}1_{A}\circ
\theta_{y}\,\mathrm{d}y\biggr)\quad\mbox{and}
\nonumber\\[-8pt]\\[-8pt]
&&P_{es}(A)=E_{ts}\biggl(\frac{1}{\alpha_{0}}\frac{1}{\overline{N}}1_{A}
\circ\eta_{0}\biggr)\qquad\mbox{for } A\in\mathcal{M}^{\infty},
\nonumber
\\
\label{eq4.5} &&\overline{N}=\frac{1}{\overline{\alpha}}
=E_{ts}\biggl(
\frac{1}{\alpha_{0}}\Big\vert\mathcal{I}\biggr)\qquad P_{es}\mbox{-a.s.},
P_{ts}\mbox{-a.s.},\mbox{ and }P\mbox{-a.s.}\\
&&\mbox{$P_{es}$ is the event-stationary PD of $P_{ts}\quad
\Leftrightarrow\quad P_{es}$ is
pseudo-ergodic.}\nonumber
\end{eqnarray}
\end{theorem}

\begin{pf}

\textit{Proof of `$\Rightarrow$'}. By Birkhoff's ergodic theorem
and the left-hand side of (\ref{eq4.3}), we obtain that, for all
$P_{es}$-integrable functions $f\dvtx M^{\infty} \rightarrow\mathbb{R}$,
the following convergence holds not only $P_{es}$-a.s. but also
$P$-a.s.:
\[
\frac{1}{n}\sum_{i=1}^{n}f\circ
\eta_{i}\rightarrow E_{es}(f| \mathcal{I}) \qquad\mbox{as }
n\rightarrow\infty.
\]
The choices $f=\alpha_{0}$ and
$f=\int_{T_{0}}^{T_{1}}1_{A}\circ\theta_{y}\,\mathrm{d}y$ (with
$A\in\mathcal{M}^{\infty}$) respectively yield that, $P$-a.s.:
\[
\frac{1}{n}T_{n}\rightarrow\overline{\alpha}\quad\mbox{and}
\quad\frac{1}{n}\int_{T(0)}^{T(n)}1_{A}
\circ\theta_{y}\,\mathrm{d}y\rightarrow E_{es}\biggl(\int
_{0}^{\alpha_{0}}1_{A}\circ\theta_{y}
\,\mathrm{d}y\Big| \mathcal{I}\biggr)\qquad\mbox{as }n\rightarrow
\infty.
\]
After replacing $n$ by $N(0, x]$ and using that $P[\overline{\alpha}
>0]=1$, we obtain that it holds $P$-a.s. that:
\begin{eqnarray*}
\frac{N(0, x]}{x}&\rightarrow&\frac{1}{\overline{\alpha}}\quad
\mbox{and}
\\
\frac{1}{N(0,
x]}\int_{0}^{x}1_{A}
\circ\theta_{y}\,\mathrm{d}y&\rightarrow& E_{es}\biggl(\int
_{0}^{\alpha_{0}}1_{A}\circ\theta_{y}
\,\mathrm{d}y\Big| \mathcal{I}\biggr)\qquad\mbox{as }x\rightarrow
\infty;
\\
\frac{N(0, x]}{x}\frac{1}{N(0, x]}\int_{0}^{x}1_{A}
\circ\theta_{y}\,\mathrm{d}y&\rightarrow&\frac{1}{\overline
{\alpha}}
E_{es}\biggl(\int_{0}^{\alpha_{0}}1_{A}
\circ\theta_{y}\,\mathrm{d}y\Big|\mathcal{I}\biggr) \qquad\mbox
{as } x
\rightarrow\infty.
\end{eqnarray*}
By taking $P$-expectations we get, again from (\ref{eq4.3}):
\[
\frac{1}{x}\int_{0}^{x}P[
\theta_{y}\varphi\in A]\,\mathrm{d}y\rightarrow P_{ts}(A):=E_{es}
\biggl(\frac{1}{\overline{\alpha}}\int_{0}^{\alpha
_{0}}1_{A}
\circ\theta_{y}\,\mathrm{d}y\biggr)\qquad\mbox{as } x\rightarrow
\infty.
\]
So, $P$ is $\operatorname{TAMS}(P_{ts})$. Note that $P=P_{es}=P_{ts}$ on
$\mathcal{I}$. Hence, $\overline{N}:=E_{ts}(N(0, 1]|\mathcal{I})$
and $1/\overline{\alpha}$ are both the
$P_{es}$-, $P$- and $P_{ts}$-a.s. limit of $N(0, x]/x$. So:
\[
\overline{N}=\frac{1}{\overline{\alpha}}\qquad P_{ts}\mbox{-a.s.}, P
\mbox{-a.s. and }P_{es}\mbox{-a.s.}\quad\mbox{and}\quad
P_{ts}[0<\overline{N}<\infty]=1.
\]

\textit{Proof of `$\Leftarrow$'}. By Birkhoff's ergodic theorem and the
right-hand side of (\ref{eq4.3}), we obtain that, for all $P_{ts}$-integrable
functions $f\dvtx M^{\infty} \rightarrow\mathbb{R}$, the following
convergence not only holds $P_{ts}$-a.s. but also $P$-a.s.:
\[
\frac{1}{x}\int_{0}^{x}f\circ
\theta_{y}\,\mathrm{d}y\rightarrow E_{ts}(f | \mathcal{I})
\qquad\mbox{as } x\rightarrow\infty.
\]
After replacing $x$ by $T_{n+1}$ and choosing
$f=1_{A}\circ\eta_{0}/\alpha_{0}$ with $A\in\mathcal{M}^{\infty
}$, we
get $P_{ts}\mbox{-a.s.}\mbox{ and} P
\mbox{-a.s.}$:
%
\begin{equation}
\label{eq4.6} \frac{1}{T_{n+1}}\int_{T_{1}}^{T_{n+1}}1_{A}
\circ\eta_{0}\circ\theta_{y}\frac{1}{\alpha_{0}\circ\theta
_{y}}\,\mathrm{d}y
\rightarrow E_{ts}\biggl(1_{A}\circ\eta_{0}
\frac{1}{\alpha_{0}} \Big| \mathcal{I}\biggr)
\qquad\mbox{as } n\rightarrow\infty.
\end{equation}
If, for $\varphi\in M^{\infty}$, $y$ is such that $T_{i}(\varphi
)\leq
y<T_{i+1}(\varphi)$, then:
\[
\alpha_{0}\circ\theta_{y}(\varphi)=\alpha_{i}(
\varphi)\quad\mbox{and}\quad1_{A}\circ\eta_{0}\circ
\theta_{y}(\varphi)=1_{A}\circ\eta_{i}(\varphi).
\]
Hence, the left-hand side of (\ref{eq4.6}) is equal to:
\[
\frac{1}{T_{n+1}}\sum_{i=1}^{n}\int
_{T_{i}}^{T_{i+1}}1_{A}\circ\eta_{i}
\frac{1}{\alpha_{i}}\,\mathrm{d}y=\frac{1}{T_{n+1}}\sum
_{i=1}^{n}1_{A}\circ\eta_{i}.
\]
Note that $N(0, x]/x\rightarrow\overline{N}$ holds $P_{ts}$- and
$P$-a.s. Replacing $x$ by $T_{n+1}$ yields that $T_{n+1}/n\rightarrow
1/\overline{N}$ also holds $P_{ts}$- and $P$-a.s. By (\ref{eq4.6}) we obtain,
$P_{ts}$- and $P$-a.s.:
\[
\frac{1}{n}\sum_{i=1}^{n}1_{A}
\circ\eta_{i}=\frac{T_{n+1}}{n}\frac{1}{T_{n+1}}\sum
_{i=1}^{n}1_{A}\circ\eta_{i}
\rightarrow\frac{1}{\overline{N}}E_{ts}\biggl(1_{A}\circ
\eta_{0}\frac{1}{\alpha_{0}} \Big| \mathcal{I}\biggr)\qquad\mbox
{as }n
\rightarrow\infty.
\]
Since $P=P_{ts}$ on $\mathcal{I}$, we obtain by taking $P$-expectation:
\[
\frac{1}{n}\sum_{i=1}^{n}P[
\eta_{i}\varphi\in A]\rightarrow P_{es}(A):=E_{ts}
\biggl(\frac{1}{\alpha_{0}}\frac{1}{\overline{N}}1_{A}\circ\eta_{0}
\biggr)\qquad\mbox{as } n\rightarrow\infty.
\]
Especially, $P=P_{es}$ on $\mathcal{I}$. For $C\in\mathcal{I}$ we have,
since $P=P_{ts}$ on $\mathcal{I}$:
\[
E(1_{C}) = E_{ts}\bigl(1_{C}/(
\alpha_{0}\overline{N})\bigr) = E_{ts}\biggl(
\frac{1}{\overline{N}}E_{ts}\biggl(\frac{1}{\alpha_{0}} \Big|
\mathcal{I}
\biggr)1_{C}\biggr) = E\biggl(\frac{1}{\overline{N}}E_{ts}\biggl(
\frac{1}{\alpha_{0}} \Big| \mathcal{I}\biggr)1_{C}\biggr).
\]
We conclude: $P$ is $\operatorname{EAMS}(P_{es})$ and
$\frac{1}{\overline{N}}E_{ts}(\frac{1}{\alpha_{0}} |
\mathcal{I})=1$ $P_{ts}$-, $P$- and $P_{es}$-a.s.
\end{pf}

\begin{remark*}
Let $P_{ts}$ be a TS point process distribution with
$\lambda_{ts}:=E_{ts}(N(0,1])<\infty$ and $P_{ts}^{0}$ the accompanying
TS Palm distribution. Note that $P_{ts}$ is $\operatorname
{TAMS}(P_{ts})$ and
$P_{ts}^{0}$ is $\operatorname{EAMS}(P_{ts}^{0})$. By Theorem \ref
{teo4.2} it follows that $P_{ts}$
is EAMS and $P_{ts}^{0}$ is TAMS. By respectively using the results
(\ref{eq4.4}), (\ref{eq2.7a}) and (\ref{eq4.5}), (\ref{eq2.7b}),
(\ref{eq2.8a}), (\ref{eq2.8c}), (\ref{eq1.1}) in 1 and 2 below, we
obtain:
\begin{enumerate}[3.]
\item[1.] $P_{ts}$ is $\operatorname{EAMS}(\tilde{P}_{ts}^{0})$ with
$\tilde{P}_{ts}^{0}(A):=\lambda_{ts}E_{ts}^{0}(\overline{\alpha}
1_{A})$ for $A\in\mathcal{M}^{\infty}$, where $\overline{\alpha}
=E_{ts}^{0}(\alpha_{0}|\mathcal{I})$;

\item[2.] $P_{ts}^{0}$ is $\operatorname{TAMS}(\tilde{P}_{ts})$ with
$\tilde{P}_{ts}(A):=E_{ts}(\overline{N}1_{A})/\lambda_{ts}$
for $A\in\mathcal{M}^{\infty}$ with
$\overline{N}=E_{ts}(N(0,1]|\mathcal{I})$;

\item[3.] $\tilde{P}_{ts}^{0}=P_{ts}^{0} \Leftrightarrow P_{ts}^{0}$
is pseudo-ergodic $\Leftrightarrow\tilde{P}_{ts}=P_{ts}$.
\end{enumerate}
\end{remark*}

The validity of EAMS (respectively, TAMS) is equivalent to the validity
of the ergodic result under the shift transformation $\eta_{1}$
(respectively under the flow $\{\theta_{y}\dvt y\in\mathbb{R}\}$).

\begin{theorem}\label{teo4.3}
Let $\Phi$ be a point process with distribution $P$
for which $P(M^{\infty} )=1$.
\begin{eqnarray*}
\mathrm{(a)}\quad&&P\mbox{ is }\mathit{EAMS} \quad\Leftrightarrow
\quad
\forall_{A\in
\mathcal{M}^{\infty} }{:}\qquad\frac{1}{n}\sum_{i=1}^{n}1_{A}
\circ\eta_{i} \mbox{ converges }P\mbox{-a.s. }(\mbox{as
}n\rightarrow
\infty);
\\
\mathrm{(b)}\quad&&P\mbox{ is EAMS with limit distribution }P_{es}
\mbox{ such that } P_{es}\bigl[0<E_{es}(\alpha_{0}
\mathcal{\vert I})<\infty\bigr]=1
\\
&&\quad\Leftrightarrow\quad\forall_{A\in\mathcal{M}^{\infty}
}{:}\qquad\frac{1}{x}N_{A}(0,x]
\mbox{ converges }P\mbox{-a.s. }(\mbox{as }x\rightarrow\infty
)\mbox{ and }
\\
&&\hspace*{33pt}\mbox{the limit of }\frac{1}{x}N(0,x]\mbox{ belongs }P
\mbox{-a.s. to }(0,\infty);
\\
\mathrm{(c)}\quad&&P\mbox{ is TAMS}\quad\Leftrightarrow\quad
\forall_{A\in
\mathcal{M}^{\infty} }{:}\qquad\frac{1}{x}\int_{0}^{x}1_{A}
\circ\theta_{y}\,\mathrm{d}y \mbox{ converges }P\mbox{-a.s.
}(\mbox{as
}x\rightarrow\infty).
\end{eqnarray*}
\end{theorem}

\begin{pf}
For (a), the implication `$\Rightarrow$' follows since
by Birkhoff's ergodic theorem the right-hand convergence holds a.s.
under the ES limit distribution of $P$, and hence under $P$ itself by
Theorem~\ref{teo4.1}. The implication `$\Leftarrow$' follows from the
theorem of
Vitali--Hahn--Saks. For `$\Rightarrow$' of (c), apply Birkhoff's ergodic
theorem to the flow $\{\theta_{y}\}$ and the TS limit distribution of
$P$, and apply Theorem~\ref{teo4.1}. For `$\Leftarrow$' of (c), apply again
Vitali--Hahn--Saks. So, only (b) is left. If $P$ is $\operatorname
{EAMS}(P_{es})$ and
$P_{es}[0<E_{es}(\alpha_{0}\mathcal{\vert I})<\infty]=1$, then
application of (a) and Theorem \ref{teo4.2} yields:\vspace*{-8pt}
%
\begin{subequation}\label{eq4.7}
\begin{eqnarray}
\label{eq4.7a}&& N_{A}(0,x]\frac{1}{N(0,x]}=\frac{1}{N(0,x]}\sum
_{i=1}^{N(0,x]}1_{A}\circ
\eta_{i}\qquad\mbox{converges }P\mbox{-a.s. }(\mbox{as
}x\rightarrow
\infty),
\\
\label{eq4.7b}&& \frac{1}{x}N_{A}(0,x]=\frac{1}{x}N(0,x]
\cdot N_{A}(0,x]\frac{1}{N(0,x]} \qquad\mbox{converges } P
\mbox{-a.s. }(\mbox{as }x\rightarrow\infty).
\end{eqnarray}
\end{subequation}
Implication `$\Leftarrow$' of (b) follows from (a), since the expression
below converges $P$-a.s.:
\[
\frac{1}{n}\sum_{i=1}^{n}1_{A}
\circ\eta_{i}=\frac{1}{n}N_{A} (0,T_{n} ]=
\frac{T_{n}}{N (0,T_{n} ]}\cdot\frac{N_{A} (0,T_{n} ]}{T_{n}}.
\]
\upqed
\end{pf}

Since $P$ is EAMS iff $P_{0}$ is EAMS, part (a) also follows from
Theorem 1
of Gray and Kieffer \cite{GK80}. By Theorem \ref{teo4.3}(a), (c) it
follows immediately
that, for point process distributions $P$ and $Q$:\vspace*{-15pt}
%
\begin{subequation}\label{eq4.8}
\begin{eqnarray}
\label{eq4.8a} Q&\ll& P\quad\mbox{and}\quad P\mbox{ is EAMS} \quad
\Rightarrow
\quad Q\mbox{ is EAMS},
\\
\label{eq4.8b} Q&\ll& P\quad\mbox{and}\quad P\mbox{ is
}\operatorname{TAMS}
\quad\Rightarrow\quad Q\mbox{ is }\operatorname{TAMS}.
\end{eqnarray}
\end{subequation}

Theorems \ref{teo4.1}--\ref{teo4.3} yield several related limit
results for AMS point
processes; we mention a few. Suppose that $P$ is $\operatorname
{EAMS}(P_{es})$ and
$P_{es}[0<E_{es}(\alpha_{0}\vert\mathcal{I})<\infty]=1$. Let $ P_{ts}$
be the TS limit-distribution of $P$ with
$\lambda_{ts}=E_{ts}(N(0,1])<\infty$ and accompanying event-stationary
PD $P_{ts}^{0}$. Then, it follows by (\ref{eq4.7a}), (\ref{eq4.7b})
that, for all
$A\in\mathcal{M}^{\infty}$:
\begin{enumerate}[(b)]
\item[(a)]
$E(\frac{1}{x}N_{A}(0,x])=\frac{1}{x}\int_{(0,x]}P^{0,y}(A)\,\mathrm
{d}\nu(y)\rightarrow\lambda_{ts}P_{ts}^{0}(A)$ as $x\rightarrow
\infty$;

\item[(b)] $E[N_{A}(0,x]/N(0,x]] \rightarrow P_{es}(A)$ while
$E(N_{A}(0,x])/E(N(0,x]) \rightarrow P_{ts}^{0}(A)$ as
$x\rightarrow\infty$; cf. (\ref{eq2.4}).
\end{enumerate}

\begin{example}\label{exe4.4}
We will use Theorem \ref{teo4.3}(a) to construct a point
process distribution $P$ which is not EAMS. Set:
\begin{eqnarray*}
a(1)&=&4\quad\mbox{and}\quad a(k)= \cases{ a(k-1), &\quad if $k$ is even,
\cr
\displaystyle\sum_{i=1}^{k-1}a(i), &\quad if
$k$ is odd} \qquad\mbox{for }k = 2, 3, \dots,
\\
b(0)&=&0\quad\mbox{and}\quad b(k)=\sum_{i=1}^{k}a(i)
\qquad\mbox{for }k = 1, 2, \dots.
\end{eqnarray*}
A sequence $(x_{i})$ of \{0, 1\}-numbers is defined as follows:
\[
x_{i}=\cases{ 1, &\quad if $i\in\bigl\{b(k)+1, \dots, b(k+1)\bigr\}$
for $k$ even,
\vspace*{1pt}\cr
0, &\quad if $i\in\bigl\{b(k)+1, \dots, b(k+1)\bigr\}$ for $k$
odd} \qquad\mbox{for }k = 0, 1, 2, \dots.
\]
Note that the sequence $(m_{n})$ with
$m_{n}=\frac{1}{n}\sum_{i=1}^{n}x_{i}$ has no limit for
$n\rightarrow\infty$ since:
\[
m_{b(2n)}\rightarrow\tfrac{1}{2}\quad\mbox{and}\quad
m_{b(2n+1)}\rightarrow\tfrac{3}{4}.
\]
A point process distribution $P$ which $P$-a.s. experiences a fixed
eventuality $A$ at the times $T_{i}$ with $x_{i}=1$ and $A^{c}$ at the
times $T_{i}$ with $x_{i}=0$, is \textit{not} EAMS.
\end{example}

Unless stated otherwise, we will always assume that the conditions about
$\overline{\alpha}$ and $\overline{N}$ in Theorem~\ref{teo4.2} are
satisfied.

\section{Absolute continuity properties equivalent to AMS}\label{sec5}

It is proven that AMS is equivalent to AC properties for $\{P_{n}\}$ and
$P^{*}$, and also to a weak AC property for $\{P^{0,x}\}$. Thus, the
class of AMS point processes is characterized in three ways. With
(\ref{eq3.12a}), (\ref{eq3.12b}), we recognize the TS subclass within
the AMS class.

If $P_{m}\ll P_{es}$ for a fixed $m\in\mathbb{Z}$, then $P_{n}\ll
P_{es}$ for all $n\in\mathbb{Z}$ since $P_{es}(A)=0$ implies
$P_{es}[\eta_{n-m}(\varphi)\in A]=0$ and hence
$P_{n}(A)=P_{m}[\eta_{n-m}(\varphi)\in A]=0$. We get, for each
$m\in\mathbb{Z}$:
\[
\{P_{n}\}\ll P_{es}\quad\Leftrightarrow_{\mathrm{def}}
\quad\forall_{n\in
\mathbb{Z} }\dvt P_{n}\ll P_{es}\quad
\Leftrightarrow\quad P_{m}\ll P_{es}.
\]

The theorem below shows that $P$ is EAMS if and only if, for one (and
hence all) $m\in\mathbb{Z}$, the intermediate distribution $P_{m}$
is absolutely continuous w.r.t. an ES point process distribution.

\begin{theorem}\label{teo5.1}
Let $P$ be a point process distribution.
Then:
\begin{enumerate}[(2)]
\item[(1)] $P \mbox{ is }\operatorname{EAMS}(P_{es})\Rightarrow\{
P_{n}\}\ll
P_{es}$,

\item[(2)] $\{P_{n}\}\ll P_{es}$ with $\delta_{-n}:=\mathrm
{d}P_{n}/\mathrm{d}P_{es} \Rightarrow P \mbox{ is }\operatorname
{EAMS}(\tilde{P}_{es})$.
\end{enumerate}
Here, $\tilde{P}_{es}(A)=E_{es}(\overline{\delta} 1_{A})$
with $\overline{\delta}:=E_{es}(\delta_{0}\vert \mathcal{I})$.
\end{theorem}

\begin{pf}
Suppose that $P$ is $\operatorname{EAMS}(P_{es})$,
and that it holds for certain $A\in\mathcal{M}^{\infty}$ and
$n\in\mathbb{Z}$ that $P_{es}(A)=0$ but $P_{n}(A)=a>0$. Set
$\tilde{A}:=A\cup(\bigcup_{k\in\mathbb{Z} }\eta_{k}^{-1}A)$. Note
that $P_{es}(\tilde{A})=0$. However,
for all $k\in\mathbb{Z}$ we have:
\[
P_{k}(\tilde{A})\geq P_{k}\bigl(\eta_{n-k}^{-1}A
\bigr)=a.
\]
Hence, $P_{es}(\tilde{A})=\lim_{m\rightarrow\infty}
\frac{1}{m}\sum_{k=1}^{m}P_{k}(\tilde{A})\geq a>0$, which
leads to a contradiction. Hence, $\{P_{n}\}\ll P_{es}$. For the reversed
implication, suppose that $\{P_{n}\}\ll P_{es}$ and
$\delta_{-n}=\mathrm{d}P_{n}/\mathrm{d}P_{es}$. For all $n\in
\mathbb{Z}$ and
$C\in\mathcal{M}^{\infty}$ it holds that
$P_{n}(C)=E_{es}(\delta_{-n}1_{C})$, but also that:
\[
P_{n}(C)=P_{-1}\bigl(\eta_{n+1}^{-1}C
\bigr)=E_{es}(\delta_{1}\cdot1_{C}\circ
\eta_{n+1})=E_{es}(\delta_{1}\circ
\eta_{-n-1}\cdot1_{C}).
\]
So, $P_{es}[\delta_{k}=\delta_{1}\circ\eta_{k-1}]=1$ for all $k\in
\mathbb{Z}$ and the discrete-time stochastic process
\{$\delta_{k}\}$ is $P_{es}$-stationary. With $\overline{\delta}
:=E_{es}(\delta_{0}\vert \mathcal{I})$ we have by
Birkhoff's ergodic theorem that, for all $ A\in\mathcal{M}^{\infty}$:
\[
\frac{1}{n}\sum_{i=1}^{n}
\delta_{-i}1_{A}\rightarrow\overline{\delta} 1_{A}
\qquad\mbox{as }n\rightarrow\infty~P_{es}\mbox{-a.s.}
\]
By taking $P_{es}$-expectation it follows that $P$ is EAMS with ES limit
distribution $\tilde{P}_{es}$.
\end{pf}

\begin{remark*}
We conclude that $P$ is AMS iff there exists an ES
point process distribution $P_{es}$ such that $P_{0}\ll P_{es}$, which
also follows by combining Theorems 2, 3 and 4 of Gray and Kieffer
\cite{GK80}. Comparison with (\ref{eq3.12a}), (\ref{eq3.12b}) learns
that the class of TS point
processes is (only) a relatively small part of the class of all AMS
point processes.
\end{remark*}

\begin{theorem}\label{teo5.2}
Let $P$ be a point process distribution, and
$P_{ts}$ a TS point process distribution with finite intensity
$\lambda_{ts}$ and associated PD $P_{ts}^{0}$. Let $P^{*}$ be defined as
in (\ref{eq3.11}). Then:
\begin{enumerate}[\textup{(b)}]
\item[\textup{(a)}] If $P\ll P_{ts}$ with $\sigma:=\mathrm
{d}P/\mathrm{d}P_{ts}$, then:
\[
P_{0}\ll P_{ts}^{0}\qquad\mbox{with }
\delta_{0}:=\mathrm{d}P_{0}/\mathrm{d}P_{ts}^{0}=
\lambda_{ts}\int_{0}^{\alpha
_{0}}\sigma\circ
\theta_{y}\,\mathrm{d}y;
\]
\item[\textup{(b)}] $P^{*}\ll P_{ts} \Leftrightarrow P_{0}\ll
P_{ts}^{0}$;

\noindent for $\sigma^{*}:=\mathrm{d}P^{*}/\mathrm{d}P_{ts}$ and
$\delta_{0}:=\mathrm{d}P_{0}/\mathrm{d}P_{ts}^{0}$ we have:
\[
P_{ts}\bigl[\sigma^{*}=\delta_{0}\circ
\eta_{0}/(\lambda_{ts}\alpha_{0})\bigr]=1\quad
\mbox{and}\quad P_{ts}^{0}\biggl[\delta_{0}=
\lambda_{ts}\int_{0}^{\alpha
_{0}}
\sigma^{*}\circ\theta_{y}\,\mathrm{d}y\biggr]=1,
\]
so $\sigma^{*}$ satisfies
$P_{ts}[\sigma^{*}\circ\eta_{0}=\sigma^{*}]=1$.
\end{enumerate}
\end{theorem}

\begin{pf}
Part (a) follows immediately by applying (\ref{eq2.8a}) and
(\ref{eq2.7a}) to $P_{ts}$, $P_{ts,0}$ and $P_{ts}^{0}$. Implication
`$\Rightarrow$' of (b) is a consequence of (a) and Theorem \ref
{teo3.1}(1). For
`$\Leftarrow$' of (b), note that by (\ref{eq3.11}), (\ref{eq2.7b})
and (\ref{eq2.8a}) we obtain
for $A\in\mathcal{M}^{\infty}$ that $P^{*}(A)$ equals:
\begin{eqnarray*}
E_{ts}^{0}\biggl(\delta_{0}/\alpha_{0}
\cdot\int_{0}^{\alpha
_{0}}1_{A}\circ
\theta_{y}\,\mathrm{d}y\biggr) &=& \frac{1}{\lambda_{ts}}E_{ts}
\biggl(\delta_{0}\circ\eta_{0}/\alpha_{0}\cdot
\frac{1}{\alpha_{0}}\int_{T(0)}^{T(1)}1_{A}\circ
\theta_{y}\,\mathrm{d}y\biggr)
\\
&=& \frac{1}{\lambda_{ts}}E_{ts}\biggl(\delta_{0}\circ
\eta_{0}\frac{1}{\alpha_{0}}1_{A}\biggr)
\end{eqnarray*}
\upqed
\end{pf}
\noindent Hence, $P$ is EAMS $\Leftrightarrow$ there exists a TS point process
distribution $P_{ts}$ such that $P^{*}\ll P_{ts}$.

\begin{theorem}\label{teo5.3}
Suppose that $P(M)=1$. For all $A\in
\mathcal{M}$, the following holds:
\begin{enumerate}[(2)]
\item[(1)] $P(A)=0 \quad\Leftrightarrow\quad P^{x}(A)=0$ for $\nu$-a.e.
$x\in\mathbb{R}$;

\item[(2)] $P_{n}(A)=0$ for all $n\in\mathbb{Z}$ with $P(F_{n})>0
\quad\Leftrightarrow\quad P^{0,x}(A)=0$ for $\nu$-a.e. $x\in\mathbb{R}$.
\end{enumerate}
\end{theorem}

\begin{pf}
The left-hand sides of (1) and (2) are, respectively,
equivalent to $\nu_{A}(B)$ being 0 for all
$B\in\operatorname{Bor}(\mathbb{R})$ and to $\mu_{A}(B)$ being 0
for all
$B\in\operatorname{Bor}(\mathbb{R})$. Next, use (\ref{eq3.1}) and
(\ref{eq3.3}).
\end{pf}

Below, we will consider the following absolute continuity properties:
\begin{eqnarray*}
&&\bigl\{P^{0,x}\bigr\}\ll P_{es}\hspace*{6pt}\quad\Leftrightarrow
_{\mathrm{def}}
\quad\mbox{for }\nu\mbox{-a.e. } x\in\mathbb{R}\mbox{ and
}\forall_{A\in \mathcal{M}}{:}
\quad P_{es}(A)=0 \quad\Rightarrow\quad P^{0,x}(A)=0
\\
&&\hspace*{6pt}\hphantom{\bigl\{P^{0,x}\bigr\}\ll P_{es}\quad
}\Leftrightarrow_{\mathrm{def}}\quad\bigl\{P^{0,x}\bigr
\}\mbox{ is \textit{absolute continuous w.r.t.} }P_{es},
\\
&&\bigl\{P^{0,x}\bigr\}\ll_{w} P_{es}\quad
\Leftrightarrow_{\mathrm{def}}\quad\forall_{A\in \mathcal
{M}}\mbox{ and for }\nu
\mbox{-a.e. } x\in\mathbb{R}{:} \quad P_{es}(A)=0\quad\Rightarrow
\quad
P^{0,x}(A)=0
\\
&&\hspace*{6.5pt}\hphantom{\bigl\{P^{0,x}\bigr\}\ll P_{es}\quad
}\Leftrightarrow_{\mathrm{def}}\quad\bigl\{P^{0,x}\bigr
\}\mbox{ is \textit{weakly absolute continuous w.r.t.} }P_{es},
\\
&&\bigl\{P^{x}\bigr\}\ll P\hspace*{19pt}\quad\Leftrightarrow
_{\mathrm{def}}\quad
\mbox{for }\nu\mbox{-a.e. } x\in\mathbb{R}\mbox{ and }\forall
_{A\in \mathcal{M}}{:}
\quad P(A)=0 \quad\Rightarrow\quad P^{x}(A)=0.
\end{eqnarray*}
The next result is an immediate consequence of Theorem \ref{teo5.3}.

\begin{corollary}\label{cor5.4}
Let $P$ be a point process distribution. It holds for event-stationary
$P_{es}$ that:
$\{P^{0,x}\}\ll P_{es} \Rightarrow \{P^{0,x}\}\ll_{w}P_{es}
\Leftrightarrow \{P_{n}\}\ll P_{es}$.
\end{corollary}

The example below shows that Theorem \ref{teo5.3} does not necessarily
imply that
$\{P^{x}\}\ll P$ and also not that $\{P^{0,x}\}\ll P_{n}$. It also shows
that $\{P^{0,x}\}\ll_{w}P_{es}$ does not necessarily imply
\mbox{$\{P^{0,x}\}\ll P_{es}$}.

\begin{example}\label{exe5.5}
Let $P$ be the distribution of an ES Poisson point
process. Note that the eventualities $A_{x}:=[\varphi\{x\}=1]$
have $P$-probability 0 and $P^{x}$-probability 1 as long as $x \neq0$.
So, $\{P^{x}\}\ll P$ is \textit{not} valid. Also $\{P^{0,x}\}\ll P_{n}$
is \textit{not} valid since $P_{n}=P$ and the eventualities
$C_{x}:=[\varphi\{-x\}=1]$ have $P^{0,x}$-probability 1 and
$P_{n}$-probability 0 for $x \neq0$. By Corollary \ref{cor5.4}, we have
$\{P^{0,x}\}\ll_{w}P$. However, $\{P^{0,x}\}\ll P$ is not valid since
for all $x \neq0$ we have $P(C_{x})=0$ while $P^{0,x}(C_{x})=1$.
\end{example}

If $P$ satisfies $P(M^{0})=1$, then $P=P_{0}$. So, the property $P\ll
P_{es}$ is not interesting for further investigation about AMS. However,
the property $P\ll P_{ts}$ is interesting since, by (\ref{eq2.8a}), it only
implies that $P^{*}\ll P_{ts}$ (i.e., no equivalence) and hence that $P$
is AMS. By Corollary \ref{cor5.4}, the property $\{P^{0,x}\}\ll
P_{es}$ also (only)
implies that $P$ is AMS. In Sections \ref{sec6}--\ref{sec8}, we will
characterize the
properties $\{P^{0,x}\}\ll P_{es}$ and $P\ll P_{ts}$ and derive
relationships between them.

\section{Absolute continuity of $\{P^{0,x}\}$ w.r.t. $P_{es}$}\label{sec6}

The property $\{P^{0,x}\}\ll P_{es}$ implies that $P$ can be expressed
in $P_{es}$. The property is stronger than AMS; we characterize it.
Below, we will use that:
\[
P_{(\eta_{n},T_{n})}(A\times B):=P[\eta_{n}\varphi\in A;T_{n}
\varphi\in B ]\qquad\mbox{for }A\in\mathcal{M}^{\infty}\mbox{
and } B\in
\operatorname{Bor}(\mathbb{R}).
\]
Suppose that $\{P^{0,x}\}\ll P_{es}$ with RN derivatives
$\{\rho_{x}\}$. Then we have, for $\nu$-a.e. $x\in\mathbb{R}$:
%
\begin{equation}
\label{eq6.1} P^{0,x}(A)=E_{es}(\rho_{x}1_{A})
\quad\mbox{for all }A\in\mathcal{M}.
\end{equation}
Since $P_{es}(M^{\infty} )=1$, it follows that $P^{0,x}(M^{\infty} )=1$
for $\nu$-a.e. $x\in\mathbb{R}$. Hence, $P(M^{\infty} )=1$ by Theorem
\ref{teo5.3}(1). By (\ref{eq3.6c}) and (\ref{eq6.1}) we can express
$P$ in $P_{es}$:
%
%
\begin{equation}
\label{eq6.2} P (A )=E_{es} \biggl(\int_{ (-T_{-k+1}, -T_{-k} ]}^{}
\rho_{y}\cdot1_{A}\circ\theta_{-y}\,\mathrm{d}\nu
(y ) \biggr),\qquad A\in\mathcal{M}^{\infty} \mbox{ and } k\in
\mathbb{Z}.
\end{equation}

\begin{theorem}\label{teo6.1}
Let $P_{es}$ be an ES distribution on
$(M^{0},\mathcal{M}^{0})$. Then:
\begin{eqnarray*}
&&\bigl\{P^{0,x}\bigr\}\ll P_{es}\mbox{ on }
\bigl(M^{0},\mathcal{M}^{0}\bigr)
\\
&&\quad\Leftrightarrow\quad\forall_{n \in\mathbb{Z}}\dvt\{
P_{(\eta_{n},T_{n})}\}\ll
P_{es}\times\nu\mbox{ on }\bigl(M^{0}\times\mathbb{R},
\mathcal{M}^{0}\otimes\operatorname{Bor}(\mathbb{R})\bigr).
\end{eqnarray*}
The RN-derivatives
\[
\rho_{x} (\varphi):=\frac{\mathrm{d}P^{0,x}}{\mathrm{d}P_{es}}
(\varphi) \quad\mbox{and}
\quad\tau_{-n} (\varphi,x ):=\frac{\mathrm{d}P_{ (\eta
_{n},T_{n} )}}{\mathrm{d}(P_{es}\times\nu)} (\varphi,x ),
\]
$\varphi\in M^{0}$, $x\in\mathbb{R}$ and $n\in\mathbb{Z}$, are
related as follows:
\begin{enumerate}[(2)]
\item[(1)] $\tau_{-n}(\varphi,x)= \rho_{x}(\varphi)\cdot1_{[T(-n)
\leq -x < T(-n+1)]}(\varphi)$ $(P_{es}\times\nu)$-a.e.

\item[(2)] $\rho_{x}(\varphi)=\sum_{k\in\mathbb{Z}}\tau
_{-k}(\varphi,x)$ $P_{es}$-a.s. for $\nu$-a.e.
$x\in\mathbb{R}$.
\end{enumerate}
\end{theorem}

\begin{pf}
The implication `$\Rightarrow$' and (1) follow
from (\ref{eq3.9}) and (\ref{eq6.1}). Next, suppose that for all
$n\in\mathbb{Z}$
the $P$-distribution of $(\eta_{n},T_{n})$ is dominated by
$P_{es}\times\nu$, with RN-derivative denoted as $\tau_{-n}(\varphi,x)$.
Set $Q^{0,x}(A):=\int_{A}(\sum_{n\in\mathbb{Z}}\tau_{-n}(\varphi
,x))\,\mathrm{d}P_{es}(\varphi)$, for $A\in\mathcal{M}^{0}$.
Note that
\[
\int_{B}Q^{0,x}\bigl(M^{0}\bigr)\,
\mathrm{d}\nu(x) = \sum_{n\in\mathbb{Z}}\int
_{B}\int_{M^{0}}\tau_{-n}(
\varphi,x)\,\mathrm{d}P_{es}(\varphi)\,\mathrm{d}\nu(x) = \nu(B)
= \int
_{B}1\,\mathrm{d}\nu(x)
\]
for all $B\in\operatorname{Bor}(\mathbb{R})$. So, for $\nu$-a.e.
$x\in\mathbb{R}$ it holds that $Q^{0,x}(M^{0})=1$ and $Q^{0,x}$ is
a probability measure. The right-hand side of (\ref{eq3.2}), with
$P^{x}(\cdot)$
replaced by $Q^{0,x}[\theta_{-x}\varphi\in\cdot]$, equals
\begin{eqnarray*}
&&\int_{\mathbb{R}}\int_{M^{0}}f(x,
\theta_{-x}\varphi)\,\mathrm{d}Q^{0,x}(\varphi)\,\mathrm{d}
\nu(x)
\\
&&\quad= \int_{\mathbb{R}}\int_{M^{0}}f(x,
\theta_{-x}\varphi)\sum_{n\in
\mathbb{Z}}
\tau_{-n}(\varphi,x)\,\mathrm{d}P_{es}(\varphi)\,\mathrm{d}
\nu(x)
\\
&&\quad= \sum_{n\in\mathbb{Z}}\int_{M^{0} \times
\mathbb{R}}f(x,
\theta_{-x}\varphi)\,\mathrm{d}P_{(\eta_{n},T_{n})}(\varphi,x)
\\
&&\quad= \sum_{n\in\mathbb{Z}}\int_{M}f
\bigl(T_{n}\varphi,\theta_{-T_{n}(\varphi
)}(\eta_{n}\varphi)
\bigr)\,\mathrm{d}P(\varphi)=\sum_{n\in\mathbb{Z}}\int
_{M}f(T_{n}\varphi,\varphi)\,\mathrm{d}P(\varphi),
\end{eqnarray*}
which is just the left-hand side of (\ref{eq3.2}). Since the family of
PD's of
$P$ is unique in the $\nu$-a.e. sense, we have $P^{0,x}=Q^{0,x}$ for
$\nu$-a.e. $x\in\mathbb{R}$. The if-part and (2) follow.
\end{pf}

\begin{corollary}\label{cor6.2}
Suppose that $\{P^{0,x}\}\ll P_{es}$ with
RN-derivatives $\{\rho_{x}\}$. Then:
\begin{enumerate}[(3)]
\item[(1)] For all $n\in\mathbb{Z}$: $P_{n}\ll P_{es}$ with
RN-derivative $\delta_{-n}=\int_{(-T(-n+1),-T(-n)]}\rho_{y}\,\mathrm
{d}\nu(y)$.

\item[(2)] For all $m\in\mathbb{Z}$ it holds that
$P_{es}[\delta_{m+1}=\delta_{m}\circ\eta_{1}]=1$, so $\{\delta
_{n}\}$ is
$P_{es}$-stationary.

\item[(3)] If it holds additionally that $P_{es}[\delta_{0}>0]=1$, then:
\[
\bigl\{P^{0,x}\bigr\}\ll P_{0}\quad\mbox{and}\quad
P_{0}\ll P_{es};\qquad\mbox{here, } \mathrm{d}P^{0,x}/
\mathrm{d}P_{0}=\rho_{x}/\delta_{0}\ P_{0}
\mbox{-a.s.}
\]
\end{enumerate}
\end{corollary}

\begin{pf}
(1) follows immediately from Theorem \ref{teo6.1}. For (2), note
that, for all $A\in\mathcal{M}^{\infty}$ and \mbox{$m\in\mathbb{Z}$}:
\[
E_{es}(1_{A}\cdot\delta_{m+1}) =
P_{-m}[\eta_{-1}\varphi\in A]=E_{es}(1_{A}
\circ\eta_{-1}\cdot\delta_{m})=E_{es}(1_{A}
\cdot\delta_{m}\circ\eta_{1}).
\]
Part (3) follows from (1) and from the fact that it holds for $\nu
$-a.e. $x\in\mathbb{R}$ that:
\[
P^{0,x}(A)=E_{es}(\rho_{x}1_{A})=E_{es}(
\rho_{x}1_{A}1_{[\delta
_{0}>0]})=E_{0}\biggl(
\rho_{x}\frac{1}{\delta_{0}}1_{A}\biggr)\qquad\mbox{for } A\in
\mathcal{M}^{0}.
\]
\upqed
\end{pf}

Note that Corollary \ref{cor6.2}(1) generalizes (\ref{eq2.7a}). By
(\ref{eq6.2}), the additional
assumption $P_{es}[\delta_{0}>0]=1$ yields that $P$ can be expressed in
terms of $P_{0}$, a property that according to (\ref{eq2.8a}) and
(\ref{eq3.11}) also
holds for TS distributions $P$ and, more generally, for distributions
$P$ with $P^{*}=P$.

\section{Absolute continuity of $P$ w.r.t. $P_{ts}$}\label{sec7}

The point processes with $P\ll P_{ts}$ are characterized within the
class of AMS point processes. We also compare the properties $P\ll
P_{ts}$, $P\ll P^{*}$, $P^{*}\ll P_{ts}$, $\nu\ll\mathrm{Leb}$, and
$P^{*}=P$. The equivalence of `$P$ is also TS' and
`$\mathcal{I}$-measurability of $\mathrm{d}P/\mathrm{d}P_{ts}$' is
proved. For
time-stationary~$P$, the property $P\ll P_{ts}$ holds equivalently for
the associated event-stationary PDs.

Assume that $P\ll P_{ts}$. Hence, $P[\varphi\{0\}=0]=1$ and
$P(M^{\infty} )=1$. Set $\sigma:=\mathrm{d}P/\mathrm{d}P_{ts}$, let
$\lambda_{ts}$ be the (finite) intensity of $P_{ts}$ and let
$P_{ts}^{0}$ be the event-stationary PD of $P_{ts}$. It follows that:
%
\begin{eqnarray}
\label{eq7.1} P(A)&=&E_{ts}(\sigma\cdot1_{A}),\qquad A\in
\mathcal{M}^{\infty},
\\
\label{eq7.2} P(A)&=&\lambda_{ts}E_{ts}^{0}\biggl(
\int_{(-T_{-k+1},-T_{-k}]}\sigma\circ\theta_{-y}
\cdot1_{A}\circ\theta_{-y}\,\mathrm{d}y\biggr),\qquad A\in
\mathcal{M}^{\infty}\mbox{ and }k\in\mathbb{Z}.
\end{eqnarray}

\begin{theorem}\label{teo7.1}
Let $P$ and $P_{ts}$ be point process distributions
and let $P^{*}$ be as in (\ref{eq3.11}); suppose that $P_{ts}$ is
time-stationary. Below, versions of RN-derivatives for $P\ll P_{ts}$,
$P^{*}\ll P_{ts}$ and \mbox{$P\ll P^{*}$} are (if existing) respectively
denoted as $\sigma$, $\sigma^{*}$ and $\tau$.
\begin{enumerate}[(b)]
\item[(a)] $P\ll P_{ts} \Leftrightarrow P\ll P^{*}$ and
$P^{*}\ll P_{ts}$; $\sigma^{*}=\frac{1}{\alpha_{0}}\int
_{T(0)}^{T(1)}\sigma\circ\theta_{y}\,\mathrm{d}y$ and $\tau=\sigma
/\sigma^{*}$;

\item[(b)] If $P\ll P_{ts}$, then: $P=P^{*} \Leftrightarrow
P_{ts}[\sigma=\sigma\circ\eta_{0}]=1$;

\item[(c)] If $P\ll P_{ts}$, then $\nu\ll\mathrm{Leb}$.
\end{enumerate}
\end{theorem}

\begin{pf}
The implication `$\Leftarrow$' of (a) is trivial. For
the implication `$\Rightarrow$', suppose that $P\ll P_{ts}$. By (\ref
{eq3.11}),
(\ref{eq7.1}), and (\ref{eq2.8a}) under $P_{ts}$, it follows that
$P^{*}\ll P_{ts}$ with
$\sigma^{*}$ as indicated. Since $\sigma^{*}\circ\eta_{0}=\sigma^{*}$,
we obtain by (\ref{eq3.11}) that not only $P^{*}[\sigma^{*}=0]=0$,
but also
$P[\sigma^{*}=0]=0$. It follows that $P\ll P^{*}$ since, because of
$P^{*}\ll P_{ts}$, we have for all $A\in\mathcal{M}^{\infty}$:
\[
E^{*}\biggl(\frac{\sigma} {\sigma^{*}}1_{A}\biggr) =
E_{ts}(\sigma1_{[\sigma
^{*}>0]}1_{A}) = P\bigl(A\cap\bigl[
\sigma^{*}>0\bigr]\bigr) = P(A).
\]
Part (b) follows from (a). For (c), suppose that $P\ll P_{ts}$. If
$B\in\operatorname{Bor}(\mathbb{R})$ satisfies $\mathrm{Leb}(B)=0$, then
$\nu(B)=0$ since:
\[
E_{ts}N(B)=\lambda_{ts}\cdot\mathrm{Leb}(B)=0\quad
\mbox{and}\quad P_{ts}\bigl[\varphi(B)=0\bigr]=1=P\bigl[\varphi(B)=0
\bigr].
\]
\upqed
\end{pf}

By Theorems \ref{teo5.2}(b) and \ref{teo7.1}(a) it follows that the
point processes with $P\ll
P_{ts}$ are just the AMS point processes for which it additionally holds
that $P\ll P^{*}$.

\subsection*{If P is time-stationary too \dots}

We will consider the consequences of $P\ll P_{ts}$ if $P$ is also
time-stationary.

\begin{theorem}\label{teo7.2}
Suppose that $P\ll P_{ts}$ and that
$\lambda_{ts}<\infty$. Then:
\begin{enumerate}[(b)]
\item[(a)] $P$ is time-stationary too $\Leftrightarrow$ there
exists an $\mathcal{I}$-measurable version of
$\mathrm{d}P/\mathrm{d}P_{ts}$.

\item[(b)] If $P$ is also time-stationary and $P_{ts}$ is ergodic, then
$P=P_{ts}$.

\item[(c)] If $P$ is also time-stationary with intensity $\lambda$ and
$P_{ts}$ is pseudo-ergodic, then $P$ is also pseudo-ergodic and
$\lambda
=\lambda_{ts}$.
\end{enumerate}
\end{theorem}

\begin{pf}
For (a), suppose that $P$ is also TS and set $\sigma
=\mathrm{d}P/\mathrm{d}P_{ts}$. By Birkhoff's ergodic theorem and
taking $E$-expectations,
we obtain for $A\in\mathcal{M}^{\infty}$:
%
\begin{eqnarray}
\label{eq7.3}&& \frac{1}{x}\int_{0}^{x}1_{A}
\circ\theta_{y}\,\mathrm{d}y\rightarrow E_{ts}(1_{A}|
\mathcal{I})\qquad\mbox{as } x\rightarrow\infty\ P\mbox{-a.s.},\\
&&P(A)=E\bigl(E_{ts}(1_{A} | \mathcal{I})
\bigr)=E_{ts}\bigl(\sigma E_{ts}(1_{A} |
\mathcal{I})\bigr)=E_{ts}(\overline{\sigma} 1_{A}).
\nonumber
\end{eqnarray}
Hence, $\overline{\sigma}:= E_{ts}(\sigma| \mathcal{I})$ is an
$\mathcal{I}$-measurable version of $\mathrm{d}P/\mathrm{d}P_{ts}$.
The if-part follows
from (\ref{eq1.1}). Statement (b) follows from (\ref{eq7.3}) since
now the limit is
$P_{ts}(A)$. For (c), note that $E(N(0,1]|\mathcal{I})$ and
$\lambda_{ts}$ are both the $P$-a.s. limit of $N(0,x]/x$ as
$x\rightarrow\infty$. Hence, $P[E(N(0,1]|\mathcal{I})=\lambda_{ts}]=1$,
$P$~is pseudo-ergodic too, and $\lambda=\lambda_{ts}$.
\end{pf}

\begin{theorem}\label{teo7.3}
Suppose that $P$ and $P_{ts}$ are both TS with
respective (finite) intensities $\lambda$ and $\lambda_{ts}$, and
accompanying event-stationary PDs $P^{0}$ and $P_{ts}^{0}$. Then:
\[
P\ll P_{ts}\quad\Leftrightarrow\quad P^{0}\ll
P_{ts}^{0}.
\]
Respective $\mathcal{I}$-measurable versions $\sigma$ and $\sigma_{0}$
of the RNs satisfy: $\lambda\sigma_{0}=\lambda_{ts}\sigma$
$P_{ts}^{0}$-a.s.
and \mbox{$P_{ts}$-a.s.}
\end{theorem}

\begin{pf}
If $P\ll P_{ts}$, then, by Theorem \ref{teo7.2}(a), we can take an
$\mathcal{I}$-measurable version $\sigma$ for $\mathrm{d}P/\mathrm
{d}P_{ts}$. By (\ref{eq2.7b}),
(\ref{eq1.1}) and (\ref{eq2.7a}) we obtain for all $A\in\mathcal
{M}^{\infty}$:
\[
P^{0}(A)=\frac{1}{\lambda} E\biggl(\frac{1}{\alpha_{0}}1_{A}
\circ\eta_{0}\biggr)=\frac{1}{\lambda} E_{ts}\biggl(\sigma
\frac{1}{\alpha_{0}}1_{A}\circ\eta_{0}\biggr)=\frac{\lambda
_{ts}}{\lambda}
E_{ts}^{0}(\sigma\cdot1_{A}).
\]
Hence, $P^{0}\ll P_{ts}^{0}$, and $\sigma_{0}=\mathrm{d}P^{0}/\mathrm
{d}P_{ts}^{0}$
satisfies $P_{ts}^{0}[\lambda\sigma_{0}=\lambda_{ts}\sigma]=1$. If
$P^{0}\ll P_{ts}^{0}$, it can be proved (as in the proof of Theorem
\ref{teo7.2}(a))
that $\sigma_{0}$ can be taken as an $\mathcal{I}$-measurable function.
By (\ref{eq2.6}), (\ref{eq2.7b}) under $P_{ts}$, and (\ref{eq2.8c})
under $P_{ts}$, we have for
$C\in\mathcal{M}^{\infty}$:
\begin{eqnarray*}
P(C)& = &\lambda E^{0}\biggl(\int_{0}^{\alpha_{0}}1_{C}
\circ\theta_{y}\,\mathrm{d}y\biggr)=\lambda E_{ts}^{0}
\biggl(\sigma_{0}\int_{0}^{\alpha
_{0}}1_{C}
\circ\theta_{y}\,\mathrm{d}y\biggr)
\\
&=&\frac{\lambda} {\lambda_{ts}}E_{ts}\biggl(\sigma_{0}\cdot
\frac{1}{\alpha_{0}}\int_{T(0)}^{T(1)}1_{C}\circ
\theta_{y}\,\mathrm{d}y\biggr) = \frac{\lambda} {\lambda_{ts}}E_{ts}(
\sigma_{0}\cdot1_{C}).
\end{eqnarray*}
\upqed
\end{pf}

\section{Relationships between absolute continuity properties}\label{sec8}

The properties $P\ll P_{ts}$ and $\{P^{0,x}\}\ll P_{ts}^{0}$ are
compared and the relationships between the accompanying RN derivatives
are investigated. If $P^{*}=P$, then $P$ is AMS iff there exist a
time-stationary $P_{ts}$ which dominates $P$.

\begin{theorem}\label{teo8.1}
Let $P$ and $P_{ts}$ be point process
distributions, where $P_{ts}$ is time-stationary, \mbox{$\lambda
_{ts}<\infty$}
and $P_{ts}^{0}$ is the accompanying Palm distribution. Then:
\[
P\ll P_{ts}\quad\Leftrightarrow\quad\nu\ll\mathrm{Leb}\quad
\mbox{and}\quad\bigl\{P^{0,x}\bigr\}\ll P_{ts}^{0}.
\]
The RN-derivatives $\sigma:=\mathrm{d}P/\mathrm{d}P_{ts}$,
$\lambda(\cdot):=\mathrm{d}\nu/\mathrm{d}\mathrm{Leb}$ and
$\rho_{x}:=\mathrm{d}P^{0,x}/\mathrm{d}P_{ts}^{0}$ satisfy:
\begin{enumerate}[(b)]
\item[(a)]
$\lambda(y)=\lambda_{ts}\cdot E_{ts}^{0}(\sigma\circ\theta_{-y})$
for Leb-a.e. $y\in\mathbb{R}$;

\item[(b)]
$P_{ts}^{0}[\lambda(y)\cdot\rho_{y}=\lambda_{ts}\cdot\sigma\circ
\theta_{-y}]=1$
for Leb-a.e. $y\in\mathbb{R}$;

\item[(c)]
$P_{ts}^{0}[\lambda(y)\cdot\rho_{y}=\lambda_{ts}\cdot\sigma\circ
\theta_{-y}$ for Leb-a.e. $y\in\mathbb{R}
]=1$;

\item[(d)] $P_{ts}[\lambda_{ts}\cdot\sigma
=\lambda(T_{k})\cdot(\rho_{T_{k}}\circ\eta_{k})]=1$ for all
$k\in\mathbb{Z}$.
\end{enumerate}
\end{theorem}

\begin{pf}
Suppose that $P\ll P_{ts}$ with RN-density $\sigma$.
First note that $\nu\ll\mathrm{Leb}$ by Theorem \ref{teo7.1}(c). Write
$\lambda(\cdot)$ for the RN-density and note that $\lambda(x)>0$ for
$\nu$-a.e. $x\in\mathbb{R}$. For
$B\in\operatorname{Bor}(\mathbb{R})$ we have:
\[
\int_{B}\lambda(x)\,\mathrm{d}x = EN(B) = \sum
_{k\in\mathbb{Z}}P[T_{k}\in B].
\]
For all $\varphi\in M^{0}$, $k\in\mathbb{Z}$, and
$y\in\mathbb{R}$ with $T_{-k}(\varphi)<-y\leq T_{-k+1}(\varphi)$,
we have $T_{k}(\theta_{-y}\varphi)=y$ and
$\eta_{k}(\theta_{-y}\varphi)=\varphi$; call this observation $(^*)$.
Taking $A=[T_{k}\in B]$ in (\ref{eq7.2}) yields
\[
\int_{B}\lambda(x)\,\mathrm{d}x = \int_{B}
\lambda_{ts}E_{ts}^{0}(\sigma\circ
\theta_{-x})\,\mathrm{d}x,
\]
which proves (a). As a consequence of (a), we have:
%
%
\begin{equation}
\label{eq8.1} \mathrm{Leb} \bigl\{x\in\mathbb{R}\dvt\lambda(x
)=0 \mbox{ and }
P_{ts}^{0} \bigl[\sigma(\theta_{-x}\varphi)\neq0
\bigr]>0 \bigr\}=0.
\end{equation}
To prove that $\{P^{0,x}\}\ll P_{ts}^{0}$, we use (\ref{eq3.2}). For
$\nu$-a.e.
$x\in\mathbb{R}$, we define probability measures $Q^{0,x}$\vspace
*{1pt} on
$(M,\mathcal{M})$ as follows:
$Q^{0,x}(C):=\lambda_{ts}E_{ts}^{0}(\sigma\circ\theta_{-x}\cdot
1_{C})/\lambda(x)$
for $C\in\mathcal{M}^{\infty}$ . By (\ref{eq7.1}), (\ref{eq2.6}),
the above observation
$(^*)$, and (\ref{eq8.1}), the left-hand side of (\ref{eq3.2}) equals:
%
\begin{eqnarray}
\label{eq8.2} E_{ts}\biggl[\sigma(\varphi)\sum
_{k\in\mathbb{Z}}f(T_{k}\varphi,\varphi)\biggr]& = & \sum
_{k\in\mathbb{Z}}\lambda_{ts}E_{ts}^{0}
\biggl[\int_{-T_{-k+1}}^{-T_{-k}}\sigma(\theta_{-x}
\varphi)f\bigl(T_{k}(\theta_{-x}\varphi),
\theta_{-x}\varphi\bigr)\,\mathrm{d}x\biggr]
\nonumber
\\
&=&\int_{-\infty}^{\infty} \int_{M}f(x,
\theta_{-x}\varphi)\lambda_{ts}\sigma(\theta_{-x}
\varphi)\,\mathrm{d}P_{ts}^{0}(\varphi)\,\mathrm{d}x
\nonumber
\\
&=&\int_{\{x\in\mathbb{R}\dvt\lambda
(x)>0\}}\int_{M}f(x,
\theta_{-x}\varphi)\lambda_{ts}\sigma(\theta_{-x}
\varphi)\frac{1}{\lambda
(x)}\,\mathrm{d}P_{ts}^{0}(\varphi)\,
\mathrm{d}\nu(x)
\nonumber
\qquad
\\
&=&\int_{-\infty}^{\infty} \int_{M}f(x,
\theta_{-x}\varphi)\,\mathrm{d}Q^{0,x}(\varphi)\,\mathrm{d}
\nu(x).
\end{eqnarray}
Note that (\ref{eq8.2}) is just the right-hand side of (\ref{eq3.2})
if we take
$P^{x}(A)=Q^{0,x}[\theta_{-x}\varphi\in A]$, $A\in\mathcal
{M}^{\infty}$. Because of the uniqueness of $\{P^{x}\}$ it follows for
$\nu$-a.e.
$x\in\mathbb{R}$ that $Q^{0,x}=P^{0,x}$, that $P^{0,x}$ is
dominated by $P_{ts}^{0}$ for $\nu$-a.e. $x\in\mathbb{R}$ and that
the RNs $\rho_{x}$ satisfy (b) with $\nu$ instead of Leb. Hence,
$\int_{\mathbb{R}}P_{ts}^{0}[\lambda(x)\rho_{x}\neq\lambda
_{ts}\sigma\circ\theta_{-x}]\,\mathrm{d}\nu(x)=0$
and
\[
\mathrm{Leb}\bigl\{x\in\mathbb{R}\dvt\lambda(x)P_{ts}^{0}
\bigl[\lambda(x)\rho_{x}\neq\lambda_{ts}\sigma\circ
\theta_{-x}\bigr]>0 \bigr\}=0.
\]
By (\ref{eq8.1}), we obtain that $P_{ts}^{0}[\lambda(x)\rho_{x}\neq
\lambda_{ts}\sigma\circ\theta_{-x}]=0$ for $\mathrm{Leb}$-a.e.
$x\in\mathbb{R}$, which proves (b). Since (b) can equivalently be
formulated as
\[
\int_{\mathbb{R}}P_{ts}^{0}\bigl[\lambda(y)
\rho_{y}\neq\lambda_{ts}\sigma\circ\theta_{-y}
\bigr]\,\mathrm{d}y=0,
\]
result (c) is just a consequence of Fubini's theorem. Next, suppose that
$\nu\ll\mathrm{Leb}$ and \mbox{$\{P^{0,x}\}\ll P_{ts}^{0}$}. By
(\ref{eq6.2}) we have,
for $A\in\mathcal{M}^{\infty}$ and $k\in\mathbb{Z}$:
%
%
\begin{equation}
\label{eq8.3} P (A )=E_{ts}^{0}\biggl(\int
_{ (-T_{-k+1}, -T_{-k} ]}^{}\rho_{y}\cdot1_{A}
\circ\theta_{-y}\lambda(y )\,\mathrm{d}y\biggr).
\end{equation}
By observation (*), we can replace $\rho_{y}$ and $\lambda(y)$ by,
respectively, $\rho_{T_{k}\circ\theta
_{-y}}\circ\eta_{k}\circ\theta_{-y}$ and
$\lambda(T_{k}\circ\theta_{-y})$. We obtain by (\ref{eq2.6}) under
$P_{ts}$ that
the right-hand side of (\ref{eq8.3}) equals:
\[
E_{ts}\bigl(\lambda(T_{k})\cdot(\rho_{T_{k}}\circ
\eta_{k})\cdot1_{A}\bigr)/\lambda_{ts}.
\]
Hence, $P\ll P_{ts}$ and (d) follows.
\end{pf}

\begin{remark*}
Theorem \ref{teo8.1} generalizes Theorem \ref{teo7.3}. Also note
that, by
(a) and (c), it holds $P_{ts}^{0}$-a.s. that:
%
%
\begin{equation}
\label{eq8.4} \rho_{x}=\frac{\sigma\circ\theta
_{-x}}{E_{ts}^{0}(\sigma\circ\theta_{-x})} \qquad\mbox{for } \nu
\mbox{-a.e. } x\in\mathbb{R}.
\end{equation}
Starting with some preliminary TS model $P_{ts}$, each measurable
function $\sigma\dvtx M^{\infty} \rightarrow[0,\infty)$ with
$E_{ts}(\sigma)=1$ can be used to transform $P_{ts}$ into a (usually)
not TS but AMS new model $P$ via
$P(A)=E_{ts}(\sigma\cdot1_{A})$, $A\in\mathcal{M}^{\infty}$ .
The accompanying family $\{P^{0,x}\}$ of shifted Palm distributions is
then dominated by the (event-stationary) Palm distribution of $P_{ts}$.
The family of RN-densities $\{\rho_{x}\}$ is given by (\ref{eq8.4}).
\end{remark*}

By Theorems \ref{teo8.1}, \ref{teo3.1}(3), and Corollary \ref
{cor5.4} it follows that for point process
distributions $P$ with $P^{*}=P$, weak absolute domination of
$\{P^{0,x}\}$ (and hence AMS) is equivalent to strong absolute
domination:

\begin{corollary}\label{cor8.2}
If $P$ is a point process distribution with
$P^{*}=P$, then:
\[
P\ll P_{ts}\quad\Leftrightarrow\quad\bigl\{P^{0,x}\bigr\}\ll
P_{ts}^{0} \quad\Leftrightarrow\quad\bigl\{P^{0,x}
\bigr\}\ll_{w}P_{ts}^{0}.
\]
\end{corollary}

The next results are immediate consequences of Corollary \ref
{cor6.2}(1), Theorem \ref{teo8.1} and~(\ref{eq3.3}).

\begin{corollary}\label{cor8.3}
Suppose that $P\ll P_{ts}$ with
$\lambda_{ts}<\infty$ and $\sigma=\mathrm{d}P/\mathrm{d}P_{ts}$. Then:
\begin{enumerate}[(c)]
\item[(a)] $\{P_{n}\}\ll P_{ts}^{0}$ with
$\delta_{-n}:=\mathrm{d}P_{n}/\mathrm{d}P_{ts}^{0}=\lambda_{ts}\int
_{T_{-n}}^{T_{-n+1}}\sigma\circ\theta_{y}\,\mathrm{d}y$ for all
$n\in\mathbb{Z}$.

\item[(b)] For $\nu$-a.e. $x\in\mathbb{R}$ and all
$A\in\mathcal{M}^{\infty}$ it holds that:
%
%
\begin{equation}
\label{eq8.5} P^{0,x} (A )=\frac{\lambda_{ts}E_{ts}^{0}(1_{A}\cdot
\sigma\circ\theta_{-x})}{\lambda_{ts}E_{ts}^{0}(\sigma\circ\theta_{-x})}.
\end{equation}
Here, numerator and denominator are just $\lambda_{A}(x)$ and
$\lambda(x)$, RN-derivatives of $\mu_{A}$ and $\nu$ with respect to Leb.
\item[(c)] If it holds additionally that
$\sigma(\varphi)=\sigma(\eta_{0}(\varphi))$ for all $\varphi\in
M^{\infty}$ (and hence $P^{*}=P$), then
$\delta_{-n}=\sigma\circ\eta_{-n}\cdot\lambda_{ts}\alpha_{-n}$.
\end{enumerate}
\end{corollary}

In the example below, results of this paper are used to investigate
consequences of a transformation via $P\ll P_{ts}$ if $P_{ts}$ is TS and
Poisson.

\begin{example}\label{exe8.4}
Let $P_{ts}$ be the distribution of a TS Poisson
point process on $\mathbb{R}$ with intensity $\lambda_{ts}$.
Suppose that $P\ll P_{ts}$ with
$\sigma(\varphi)=\lambda_{ts}\alpha_{0}(\varphi)/2$ for $\varphi
\in
M^{\infty}$. It follows that, for $x\geq0$:
\[
P[\alpha_{0}>x]=\mathrm{e}^{-\lambda
_{ts}x}\bigl(\lambda_{ts}^{2}x^{2}/2+
\lambda_{ts}x+1\bigr)\quad\mbox{and}\quad P_{ts}[
\alpha_{0}>x]=\mathrm{e}^{-\lambda_{ts}x}(\lambda_{ts}x+1).
\]
Hence, $\alpha_{0}$ is under $P$ stochastically larger than under
$P_{ts}$. However, for $i\neq0$ the distributions of $\alpha_{i}$ under
$P$ and $P_{ts}$ are the same. By Theorems \ref{teo7.1}(b) and \ref
{teo3.1}(2) it follows
that, as under~$P_{ts}$, it holds under $P$ that $T_{1}$ is
(conditionally) uniformly distributed on $(0,\alpha_{0})$. Starting with
Theorem~\ref{teo8.1}(a), (\ref{eq8.4}) and (\ref{eq8.5}), we obtain
after tough calculations that,
for $x\in\mathbb{R}$:
\begin{eqnarray*}
&&\lambda(x) =\lambda_{ts}^{2}E_{ts}^{0}(
\alpha_{0}\circ\theta_{-x})/2 = \lambda_{ts}-
\lambda_{ts}\exp(-\lambda_{ts}|x|)/2,
\\
&&\rho_{x}=\lambda_{ts}\alpha_{0}\circ
\theta_{-x}/\bigl(2-\exp\bigl(-\lambda_{ts}|x|\bigr)\bigr),
\\
&&P^{0,x}(A)=\lambda_{ts}E_{ts}^{0}(1_{A}
\cdot\alpha_{0}\circ\theta_{-x})/\bigl(2-\exp\bigl(-
\lambda_{ts}|x|\bigr)\bigr).
\end{eqnarray*}
The independence of the interval lengths under $P_{ts}^{0}$ yields that:
\begin{itemize}
\item[-] if $x\leq0$ and $A\in\sigma\{\alpha_{-1},\alpha
_{-2},\dots\}$
then $P^{0,x}(A)=P_{ts}^{0}(A)$;

\item[-] if $x>0$ and $A\in\sigma\{\alpha_{0},\alpha_{1},\dots\}$ then
$P^{0,x}(A)=P_{ts}^{0}(A)$.
\end{itemize}

More generally, it can be proven that, for each $n\in\mathbb{Z}$,
the RN-derivative $\sigma_{n}:=\alpha_{n}/E_{ts}(\alpha_{n})$ transforms
the TS Poisson distribution $P_{ts}$ into an AMS distribution $P$ that
conserves the independence of the interval lengths $\alpha_{k}$ and for
$k\neq n$ also their distributions, but making $\alpha_{n}$
stochastically larger. However, an RN-derivative of the form $\sigma
=\gamma_{0}\alpha_{0}+\gamma_{1}\alpha_{1}$ with $\gamma_{0}\geq0$,
$\gamma_{1}\geq0$ and $\gamma_{1}=\lambda_{ts}-2\gamma_{0}$ transforms
$P_{ts}$ into a distribution under which $\alpha_{0}$ and $\alpha_{1}$
are independent if and only if
$(\gamma_{0},\gamma_{1})=(0,\lambda_{ts})$ or
$(\gamma_{0},\gamma_{1})=(\lambda_{ts}/2, 0)$.
\end{example}

\section*{Acknowledgement}

The author would like to thank the referees for their helpful and
detailed comments.


%

\printhistory

\end{document}